\documentclass[%
 aip,
 amsmath,amssymb,
 reprint
]{revtex4-1}

\usepackage{graphicx}
\usepackage{dcolumn}
\usepackage{bm}
\usepackage[utf8]{inputenc}
\usepackage[T1]{fontenc}
\usepackage{mathptmx}
\usepackage{etoolbox}

\makeatletter
\def\@email#1#2{%
 \endgroup
 \patchcmd{\titleblock@produce}
  {\frontmatter@RRAPformat}
  {\frontmatter@RRAPformat{\produce@RRAP{*#1\href{mailto:#2}{#2}}}\frontmatter@RRAPformat}
  {}{}
}%
\makeatother
                                            
\usepackage{mathtools}
\usepackage{physics}

\usepackage{csquotes}   
\MakeOuterQuote{"}    

\DeclarePairedDelimiterX\innerp[2]{\langle}{\rangle}{#1,#2}
\DeclarePairedDelimiterX\pair[2]{\langle}{\rangle}{#1,#2}
\DeclarePairedDelimiterX\poisson[2]{\{}{\}}{#1,#2}

\newcommand{\Gtotal}{\Gamma_{\rm tot}}

\newcommand{\vv}{\vb{v}}
\newcommand{\ii}{\mathrm{i}}

\newcommand{\Ymin}{Y_{\mathrm{min}}}

\newcommand{\cR}{\mathcal{R}}
\newcommand{\bcR}{\bar{\mathcal{R}}}

\renewcommand{\r}{\mathbf{r}}
\newcommand{\R}{\mathbf{R}}

\DeclareMathOperator{\sign}{sign}

\newcommand{\Equilibrium}{\mathcal{E}}
\newcommand{\Singularity}{\mathcal{S}}
\newcommand{\EqTriPM}{\Equilibrium_{\rm tri}^\pm}

\newcommand{\EqMinus}{\Equilibrium_{-1}}
\newcommand{\EqPlus}{\Equilibrium_{1}}
\newcommand{\EqGamma}{\Equilibrium_{\Gamma}}
\newcommand{\SingPP}{\Singularity_{11}}
\newcommand{\SingPGamma}{\Singularity_{1\Gamma}}
\newcommand{\SingMGamma}{\Singularity_{-1\Gamma}}
\newcommand{\DomainPhysical}{{\mathcal{D}_\mathrm{phys}}}
\newcommand{\rhocritminus}{\rho^-_{\mathrm{crit}}}
\newcommand{\rhocritplus}{\rho^+_{\mathrm{crit}}}

\graphicspath{{Images/}}

\begin{document}

\preprint{AIP/123-QED}

\title{A new canonical reduction of three-vortex motion and its application to vortex-dipole scattering}

\author{A. Anurag}

\author{R. H. Goodman}%
 \email{goodman@njit.edu}

\author{E. K. O'Grady}

\affiliation{Department of Mathematical Sciences, 
             New Jersey Institute of Technology,
             323 Martin Luther King Blvd.,
             Newark NJ, 07102}%
\date{\today}

\begin{abstract}
We introduce a new reduction of the motion of three point vortices in a two-dimensional ideal fluid. This proceeds in two stages: a change of variables to  Jacobi coordinates and then a Nambu reduction. The new coordinates demonstrate that the dynamics evolve on a two-dimensional manifold whose topology depends on the sign of a parameter $\kappa_2$ that arises in the reduction. For $\kappa_2>0$, the phase space is spherical, while for $\kappa_2<0$, the dynamics are confined to the upper sheet of a two-sheeted hyperboloid. We contrast this reduction with earlier reduced systems derived by Gröbli, Aref, and others in which the dynamics are determined from the pairwise distances between the vortices. The new coordinate system overcomes two related shortcomings of Gröbli's reduction that have made understanding the dynamics difficult: their lack of a standard phase plane and their singularity at all configurations in which the vortices are collinear. We apply this to two canonical problems. We first discuss the dynamics of three identical vortices and then consider the scattering of a propagating dipole by a stationary vortex. We show that the points dividing direct and exchange scattering solutions correspond to the locations of the invariant manifolds of equilibria of the reduced equations and relate changes in the scattering diagram as the circulation of one vortex is varied to bifurcations of these equilibria.
\end{abstract}

\maketitle

\section{Introduction}
\label{sec:intro}
The mutually induced motion of point vortices in a two-dimensional inviscid incompressible fluid is a classical topic in fluid mechanics. The positions of the vortices are described by a Hamiltonian system of ordinary differential equations that has been well-studied for over 150 years\cite{Newton.2001,Meleshko:2007}. These ODEs remain relevant because of their deep connection to turbulence in Bose-Einstein condensates and other quantum fluids, as summarized, for example, by Lydon\cite{Lydon.2022}.

The point-vortex model idealizes a near-two-dimensional inviscid incompressible fluid in which the vorticity is confined to a finite number of discrete points. Each such point vortices induces a velocity that, in turn, causes the vortices to move. It is a standard topic in elementary fluid mechanics textbooks\cite{Chorin:1993} and is well covered in Newton's textbook devoted to the subject\cite{Newton.2001}. 

Systems of three vortices are the smallest systems with time-dependent inter-vortex distances. They are integrable yet display various behaviors depending on the three circulations. Solutions to system~\eqref{N_vortex_eqns} evolve in a $2N$-dimensional phase space, so reducing the dimensionality is necessary to understand the dynamics. This paper aims to introduce a geometric reduction to the three-vortex problem that avoids introducing artificial singularities in the dynamics. Previously used reductions introduce such singularities because they are incompatible with the topology of the manifold on which the dynamics occur. This has made reasoning about the dynamics more difficult because the singularities get in the way of applying standard geometric phase-space arguments. 

We apply this reduction to two cases of the three-vortex problem: the motion of three identical vortices and the scattering of a propagating dipole by a third, initially stationary vortex. In each case, the new form of the equations dramatically simplifies the application of dynamical systems reasoning. 

Helmholtz derived the model of point-vortex motion describing the interaction of $N$ point vortices, defined by the system of $2N$ ODEs in 1858\cite{Helmholtz:1858}:
\begin{equation} \label{N_vortex_eqns}
   \dv{x_i}{t}  = - \sum_{j \neq i}^N \Gamma_j \frac{\left(y_i-y_j\right)}{\norm{\r_i-\r_j}^2}, \quad 
   \dv{y_i}{t}  =   \sum_{j \neq i}^N \Gamma_j \frac{\left(x_i-x_j\right)}{\norm{\r_i-\r_j}^2}.
\end{equation}
Here $\r_i = \pair{x_i}{y_i}$ denotes the position of the $i$th point vortex and $2\pi\Gamma_i$ represents its circulation.
The equations conserve an energy
\begin{equation} \label{N_vortex_hamiltonian}
    H(\r_1,\ldots,\r_N)=-\frac{1}{2} \sum_{i< j}^N  \Gamma_i \Gamma_j   \log \norm{ \r_i-\r_j }^2.
\end{equation} 
In 1876, Kirchhoff noted that system~\eqref{N_vortex_eqns} has a Hamiltonian formulation\cite{Kirchhoff:1876},
\begin{equation} \label{N_vortex_ham_eqns}
  \dv{x_i}{t}  = \frac{1}{\Gamma_i}\pdv{H}{y_i}, \quad 
  \dv{y_i}{t}  = -\frac{1}{\Gamma_i}\pdv{H}{x_i}.
\end{equation}

System~\eqref{N_vortex_eqns} has three well-known conservation laws, which we write as 
\begin{equation}\label{constants_of_motion}
    \vb{M}=\pair{M_x}{M_y}= \sum_{i=1}^N\Gamma_i \r_i,  \qand
    \Theta =\sum_{i=1}^N \Gamma_i \norm{\r_i}^2.
\end{equation}
The quantities $\vb{M}$ and $\Theta$ are called the \emph{linear impulse} and the \emph{angular impulse}, respectively. In the case that $\Gtotal = \sum_{i=1}^N \Gamma_i \neq 0$, then 
\begin{equation}\label{center_of_vorticity}
\r_0 = \vb{M}/\Gtotal
\end{equation}
defines the location of the conserved center of vorticity. In this case, taking $\r_0$ at the origin is natural.

The paper is organized as follows. 
In Sec.~\ref{sec:grobli}, we introduce the reduced equations that Gröbli used to integrate the equations of motion and Aref's interpretation of this system as trilinear coordinates for $\mathbb{R}^2$. This section concludes by discussing other reductions of the three-vortex system in the existing literature.
In Sec.~\ref{sec:Nvortex}, we review some ideas from Hamiltonian mechanics, including the Poisson bracket.
Sec.~\ref{sec:reduction} describes the reduction techniques, introducing with Jacobi coordinates in Sec.~\ref{sec:jacobi} and Nambu brackets in Sec.~\ref{sec:nambu} before applying these two methods to the three-vortex system in Sec.~\ref{sec:reduce3vortex}.
After this, we apply the reduced system to explore three cases of vortex motion. First, in Sec.~\ref{sec:111}, we consider the canonical case of three identical vortices.
Sec.~\ref{sec:11m1}, considers vortices with circulation $(1,1,-1)$ in which a vortex dipole is scattered by a third, initially stationary vortex of the same absolute circulation. This section contains a review of the scattering problem. Using the reduced equations, we derive an evolution equation for the instantaneous scattering angle. This is integrated in Appendix~\ref{sec:appendix}. We explain the scattering behavior, including the critical transition between direct and exchange scattering, entirely in terms of phase planes of the reduced problems.
In Sec.~\ref{sec:generalization}, we extend the analysis to the case where the initially stationary vortex has circulation $\Gamma \neq 1$.
We conclude in Sec.~\ref{sec:conclusion} with a discussion of the possible future applications of the coordinate reduction method.

\section{Gröbli's reduction and trilinear coordinates}
\label{sec:grobli}
Gröbli's 1877 doctoral thesis was the first to explore the complex dynamics that can arise in systems of three or more vortices. He simplified the three-vortex problem by deriving evolution equations for the pairwise distances between vortices\cite{Grobli:1877}, finding that these satisfy
\begin{equation}\label{l_eqn}
\dv{t} \begin{pmatrix} \ell_{23}^2 \\ \ell_{31}^2 \\ \ell_{12}^2 \end{pmatrix} =
4 \sigma   A 
\begin{pmatrix} 
   \Gamma_1 \left(\ell_{12}^{-2} -\ell_{31}^{-2} \right); \\
   \Gamma_2 \left(\ell_{23}^{-2} -\ell_{12}^{-2} \right);\\
   \Gamma_3 \left(\ell_{31}^{-2} -\ell_{23}^{-2} \right),
\end{pmatrix}\end{equation}
where $\ell_{ij}$ is the distance between vortices $i$ and $j$, $A$ is the area of the triangle formed by the vortices, and $\sigma=\pm 1$ gives the orientation of the triangle spanned by the three vortices visited in numerical order, taking the value $+1$ if they appear in clockwise order and $-1$ if counterclockwise. Since the area of a triangle can be obtained from the side lengths using Heron's formula, 
\begin{equation*}
A = \frac{1}{4}{\left(2 \ell_{12}^2 \ell_{23}^2 +2 \ell_{23}^2 \ell_{31}^2 + 2 \ell_{31}^2 \ell_{12}^2 - \ell_{12}^4 - \ell_{23}^4 - \ell_{31}^4\right)}^{1/2},
\end{equation*}
this is a closed system. 

System~\eqref{l_eqn} leads easily to a conservation law 
\begin{equation}\label{grobli_conservation}
\Gamma_1 \Gamma_2 \ell_{12}^2 + \Gamma_2\Gamma_3 \ell_{23}^2 + \Gamma_3\Gamma_1 \ell_{31}^2 = 3 L \Gamma_1 \Gamma_2\Gamma_3.
\end{equation}
The constant $L$ may also be obtained by an appropriate combination of the constants defined in Eq.~\eqref{constants_of_motion} and is proportional to $\Theta$ if the center of vorticity $\r_0$ is taken at the origin. Depending on the strengths of the three vortices, this quadratic invariant may or not be positive definite, which has consequences for the dynamics. Using this conservation law to eliminate one variable, say $\ell_{31}$, Gröbli reduced the system to quadratures. This system and others derived from it are used in most subsequent studies of the three-vortex problem \cite{Aref.1979, Synge.1949, Novikov.1975, Tavantzis.1988, Blackmore.2007, Muller.2014, Kallyadan:2022, Reinaud:2021,Reinaud:2022, Gotoda:2021, Hirakui:2021}, a history that Aref and his collaborators researched extensively\cite{Aref.1992, Meleshko:2007}. We have translated Gröbli's dissertation into English and posted it on arXiv.org\cite{Goodman:2024}.

\begin{figure*}[t] 
   \centering
   \includegraphics[height=3in]{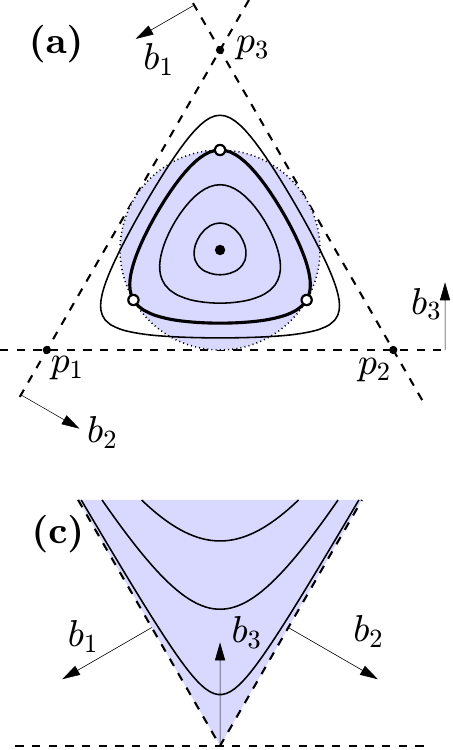} 
   \includegraphics[height=3in]{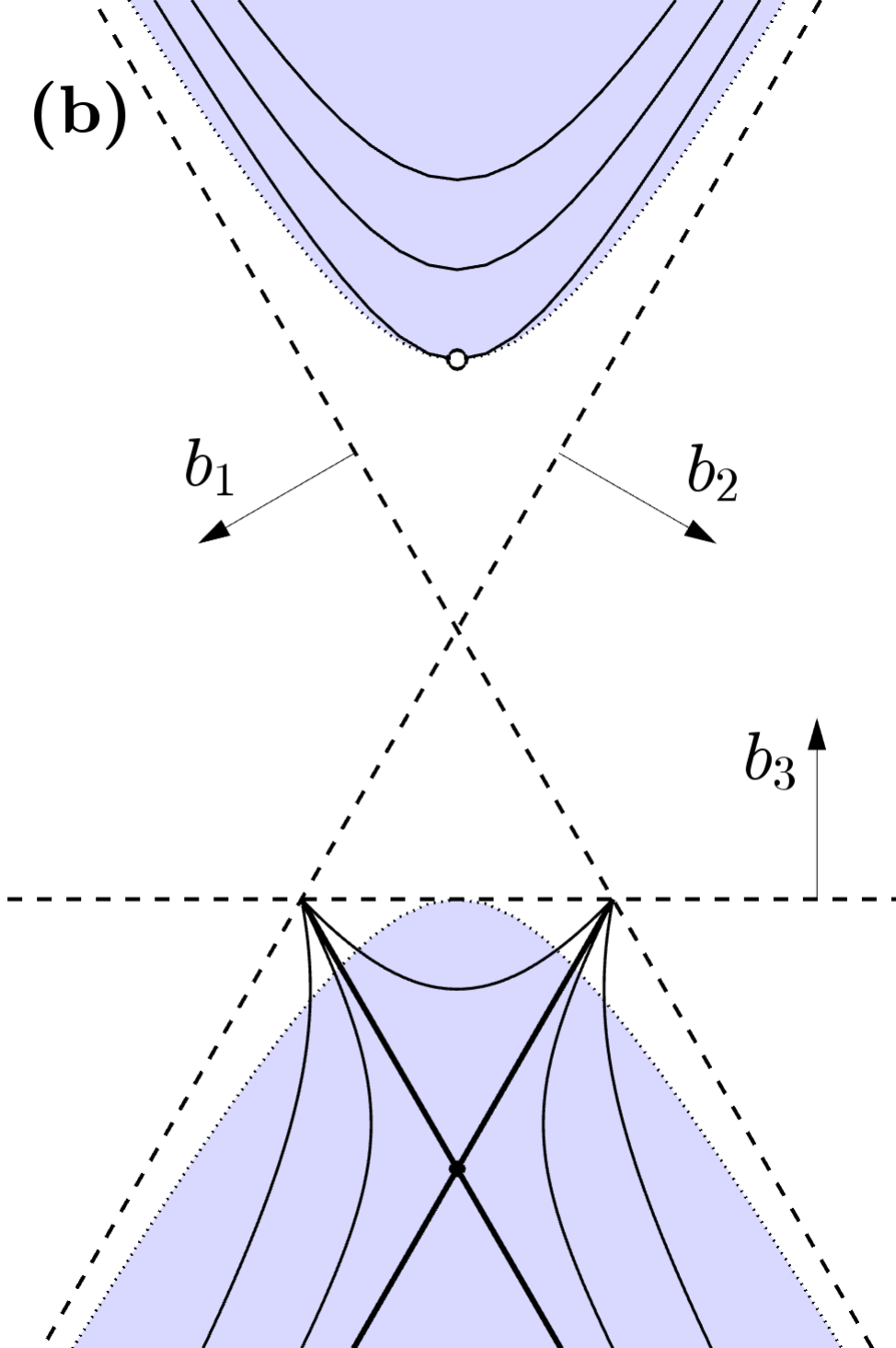} 
   \caption{Phase diagrams in trilinear coordinates. Solid black lines are level sets of the Hamiltonian, but only the portions of these curves in $\DomainPhysical$, shown here as shaded regions, are meaningful. Portions of the level sets lying outside them have no physical meaning. Collinear relative equilibria are marked $\circ$ and equilateral triangle relative equilibria $\bullet$. Heavier curves are separatrices. \textbf{(a)} Vortices of circulations $(1,1,1)$. \textbf{(b)} Vortices of circulation $(1,1,-1)$ with $L\neq0$. \textbf{(c)} Vortices of circulation $(1,1,-1)$ with $L=0$.}
   \label{fig:trilinear_diagram}
\end{figure*}

For $L\neq0$, Aref defines new variables
\begin{equation}\label{bj}
b_1=\frac{l_{23}^2}{\Gamma_1 L}, \quad b_2=\frac{l_{31}^2}{\Gamma_2 L}, \quad b_3=\frac{l_{12}^2}{\Gamma_3 L},
\end{equation}
which must then satisfy, by Eq.~\eqref{grobli_conservation},
$$
b_1+b_2+b_3 = 3.
$$
These may be interpreted as trilinear coordinates for the plane. That is, given three points $p_1$, $p_2$, and $p_3$ that form an equilateral triangle of height 3, any point in the plane is uniquely specified by the triplet of signed distances  $b_j$ from this point to the lines containing sides $j$ of the triangle, as illustrated in Fig.~\ref{fig:trilinear_diagram}(a). This has precedent in earlier works of Synge and Novikov\cite{Novikov.1975,Synge.1949}, which use a trilinear coordinate system somewhat different from Aref's.

For $L=0$, we may omit the factor of $L^{-1}$ from the definition of the trilinear coordinates in Eq.~\eqref{bj} and find instead
$$
b_1+b_2+b_3 = 0.
$$

The dynamics of the trilinear coordinates $b_j$ describe the motion of a point in the plane. Since coordinates $\ell_{ij}$ represent the sides of a triangle, they must satisfy the triangle inequality, and not all triples represent physical configurations that satisfy this constraint.  Let $\DomainPhysical \subset \mathbb{R}^3 $ denote the domain of physical configurations. Each point on its interior represents two distinct phase points related by mirror symmetry.  The boundary $\partial \DomainPhysical$ consists of collinear configurations of the three vortices, and Aref showed it describes a conic section in the plane. It is an ellipse for certain sets of circulations; for others, it is a hyperbola. Fig.~\ref{fig:trilinear_diagram} shows three such images. Panel~(a) is the phase diagram for circulation values $(1,1,1)$, the trajectories (level sets of a rescaled Hamiltonian denoted by $\theta$) are confined to lie inside the circle, which is interior and tangent to the triangle formed by the three axes of the trilinear coordinate system. The phase diagram for circulations $(1,1,-1)$ is shown in the phase diagram for circulations $(1,1,-1)$ for $L\neq0$ in panel~(b) and for $L=0$ in panel~(c). In the first, $\DomainPhysical$ is bounded by a circle; in the second, a hyperbola; and in the third, the $b_1$ and $b_2$ axes. These correspond to Figures 2, 4, and 5 in Ref.~\onlinecite{Aref.1979}.

We briefly summarize the reasoning Aref uses to interpret these figures. The evolving vortex configurations move along level sets of the Hamiltonian in the shaded $\DomainPhysical$ regions. Components of these level sets outside $\DomainPhysical$ have no physical meaning. When a trajectory reaches $\partial\DomainPhysical$, the three vortices are collinear. The motion passes through the collinear configuration and continues on the phase diagram by reversing direction and retracing its path. Thus, each point in $\DomainPhysical$ corresponds to two configurations of opposite orientation. Points where trajectories are tangent to $\partial\DomainPhysical$ correspond to collinear relative equilibria, and the orbits connected to them are their stable and unstable manifolds.

Panel~(b) shows the phase diagram for circulations $(1,1,-1)$ and $L\neq0$, for which $\partial \DomainPhysical$ is a hyperbola. The dotted portions of curves lying outside the boundary are non-physical. Trajectories that cross $\partial \DomainPhysical$ transversely immediately reverse direction and retrace the same path. Trajectories tangent to $\partial \DomainPhysical$ from the interior are invariant manifolds, and their points of tangency are hyperbolic fixed points. When $L=0$, the triangle in panel (b) shrinks to a point, and the two regions bounded by hyperbolas become wedges. The dynamics on the upper wedge of $\DomainPhysical$ is shown in panel~(c).

Points on $\partial \DomainPhysical$ in Fig.~\ref{fig:trilinear_diagram} correspond to collinear arrangements. 
Such arrangements are common: many families of periodic orbits pass through such states twice per period, and three of the five possible rigidly rotating configurations of three vortices are collinear. The evolution equations are singular on $\partial \DomainPhysical$ due to the square root that appears in Heron's formula.
Thus, linearization fails, and even finding the linear stability of the collinear states is difficult. The singularity of the reduced ODE system is an artifact of the reduced coordinate system. It is not present in the vortex motion equations~\eqref{N_vortex_eqns}, which are singular only at singularities of the Hamiltonian~\eqref{N_vortex_hamiltonian}, i.e., when two or more vortices occupy the same location. 

The images in Figure~\ref{fig:trilinear_diagram} are phase planes. Still, the above considerations show that reading the dynamics from this phase plane takes more effort than from a standard one. Certain information, like the stability of collinear fixed points, is not obtainable in this representation. Previous studies have approached different aspects of three-vortex dynamics using various reduction approaches. The first is Conte's 1979 \emph{Thèse d'État}, which appeared only as a technical report until its 2015 publication\cite{conte2015exact}. This reduces the system to an evolution equation for a single complex parameter $\zeta$, which, according to the authors, describes the shape of the triangle formed by the three vortices. The authors study many aspects of the dynamics using this reduction, but we have found the change of variables difficult to interpret. Tavantzis and Ting used Synge's trilinear coordinates to make a detailed study of the dynamics dependence on the circulation of the three vortices\cite{Tavantzis.1988}, but this has similar problems to Aref's trilinear formulation. Aref returned to the stability of collinear arrangements and introduced yet another reduction\cite{Aref.2009c5n}. This algebraic approach describes only the relative equilibria and does not apply to the dynamics more broadly.

Other reductions have included the \emph{angles} between the triangle's edges. For example, Krishnamurthy derived a system for the three angles plus the radius of the circle circumscribing the triangle\cite{Krishnamurthy:2018}, and Makarov derived a phase-plane representation for one side-length and one angle\cite{Makarov.2021}. Stremler derived an especially useful system of equations, noting that since the interior angles of a triangle must sum to $2 \pi$, they can be used as a trilinear coordinate system\cite{Stremler:2021}. Since all interior angles must be positive, the physical domain coincides with the triangle's interior, and collinear states occur at the triangle's vertices. The vertices are singular since all collinear configurations with the same central vortex degenerate to a single point, including any collinear relative equilibria. Thus, this coordinate system runs into difficulties near collinear arrangements, mirroring the weakness of trilinear coordinates.

\section{Further Mathematical Preliminaries}
\label{sec:Nvortex}

The Poisson bracket used to describe the dynamics of $N$ point vortices is defined by
\begin{equation} \label{poisson_bracket}
\poisson{F(\r)}{G(\r)} = \sum_{i=1}^N \frac{1}{\Gamma_i} \left( \pdv{F}{x_i} \pdv{G}{y_i} - \pdv{F}{y_i} \pdv{G}{x_i} \right),
\end{equation}
where 
$$
\r = \binom{\mathbf{x}}{\mathbf{y}}, \qand \mathbf{x}, \mathbf{y} \in \mathbb{R}^N.
$$
Together with the chain rule, this implies that if $\r$ evolves according to Eq.~\eqref{N_vortex_ham_eqns}, any function $F(\r)$  evolves according to
\begin{equation}\label{Fdot}
\dv{F}{t} = \poisson{F}{H}.
\end{equation}

Due to the factor of $\frac{1}{\Gamma_i}$ in these equations,  Hamiltonian system~\eqref{N_vortex_ham_eqns} is not in canonical form. It may be canonically normalized by introducing variables 
\begin{equation}\label{xy_to_pq}
    q_i=\sqrt{\abs{\Gamma_i }} x_i \qand 
    p_i=\sqrt{\abs{\Gamma_i }} \sign(\Gamma_i) y_i,
\end{equation}
which renders both the equations and the Poisson brackets into the standard forms
\begin{equation}\label{dqdt}
    \dv{q_i}{t} = \pdv{H}{p_i}, \quad 
    \dv{p_i}{t} = -\pdv{H}{q_i}.
\end{equation}
and 
\begin{equation*}
\poisson{F}{G} = \sum_{i=1}^N  \left( \pdv{F}{q_i} \pdv{G}{p_i} - \pdv{F}{p_i} \pdv{G}{q_i} \right).
\end{equation*}

\subsection{\label{sec:3vortex}The three-vortex system}

In this paper, we confine our attention to a system of three vortices with non-vanishing total circulation. The dynamics conserve the Hamiltonian
\begin{multline*}
H = -\frac{\Gamma_1 \Gamma_2}{2} \log{\norm{\r_1 - \r_2}}^2 
    - \frac{\Gamma_1 \Gamma_3}{2} \log{\norm{\r_1 -\r_3}}^2 \\
    - \frac{\Gamma_2 \Gamma_3}{2} \log{\norm{\r_2 -\r_3}}^2,
\end{multline*}
 the center of vorticity
\begin{equation*} 
\r_0 = \frac{\Gamma_1 \r_1 + \Gamma_2 \r_2 + \Gamma_3 \r_3}
             {\Gamma_1+ \Gamma_2 + \Gamma_3},
\end{equation*}
and the angular impulse
\begin{equation*} 
\Theta =  \Gamma_1 \norm{\r_1}^2 + \Gamma_2 \norm{\r_2}^2 + \Gamma_3 \norm{\r_3}^2.
\end{equation*}
Without loss of generality, we can assume that 
\begin{equation}
\label{Gamma_signs}
\Gamma_1 \ge \Gamma_2 \ge \Gamma_3 \qand \Gamma_2 >0.
\end{equation}

\section{ Reduction by stages} 
\label{sec:reduction}

Choosing the proper coordinate system may significantly clarify the study of a particular phenomenon, but how to construct such a coordinate system may not be obvious. In the ideal case, any such coordinates should be easy to interpret, which implies, among other things, that they should have a clear meaning and be invertible so that we may reconstruct the original motion from the transformed motion.  We insist on using \emph{canonical} changes of variables, those which preserve the Hamiltonian form of the evolution equations. The reduction proceeds in two steps. The first change of variables to Jacobi coordinates is canonical, while the second change, a Nambu reduction, requires us to generalize the framework of Hamiltonian systems. In between, we apply a normalization of the form~\eqref{xy_to_pq}.

\subsection{\label{sec:jacobi}Jacobi coordinate reductions}
Jacobi coordinates are a standard tool for reducing $n$-body problems, especially in celestial mechanics, and are discussed at length in Jacobi's 1866 Lectures on Dynamics\cite{Jacobi1866}. The method is straightforward and underlies the reductions used in many studies of vortex interactions\cite{Smith.2011}. Still, the only point-vortex paper we have found that references the method by name is a recent one by Luo et al.\cite{Luo:2022}.

The Jacobi coordinate transformation consists of iteratively applying the change of variables
\begin{equation} \label{Jacobi1}
\begin{aligned}
 \tilde{\r}_1 &= \r_1 - \r_2; & \tilde{\Gamma}_1 & = \frac{\Gamma_1 \Gamma_2}{\Gamma_1 + \Gamma_2}; \\
 \tilde{\r}_2 & = \frac{\Gamma_1 \r_1 + \Gamma_2 \r_2}{\Gamma_1+\Gamma_2}; & \tilde{\Gamma}_2 & = \Gamma_1 + \Gamma_2,
\end{aligned}
\end{equation}
where $\Gamma_1 + \Gamma_2 \neq 0$.

The variables $\tilde{\Gamma}_1$ and $\tilde{\Gamma}_2$, are known, respectively, as the \emph{reduced circulation} and \emph{total circulation} of the pair.

Solving for $\r_1$ and $\r_2$ and substituting these values back into the Hamiltonian~\eqref{N_vortex_hamiltonian} yields a new Hamiltonian $H(\tilde{\r}_1,\tilde{\r}_2,\r_3,\ldots,\r_N)$ and evolution equations of the form~\eqref{N_vortex_ham_eqns} with circulations $\tilde{\Gamma}_1,\tilde{\Gamma}_2,\Gamma_3,\ldots,\Gamma_N$. We then apply a similar transform to Eq.~\eqref{Jacobi1} to $\tilde{\r}_2$, $\r_3$, $\tilde{\Gamma}_2$ and $\Gamma_3$, repeating the process for each pair until $\r_N$ has been transformed. The transformed circulations redefine the Poisson bracket~\eqref{poisson_bracket} and thus the evolution equations~\eqref{Fdot}.

To prevent division by zero in the reduction procedure, we assume that 
\begin{equation} \label{partial_sum_nonzero}
\sum_{j=1}^k \Gamma_j \neq 0, \text{ for all } k \le N.
\end{equation}
For $k=N$, this represents an assumption about the set of vortices, while for $k<N$, it is merely an assumption about their labels' ordering and is consistent with assumption~\eqref{Gamma_signs} and~\eqref{Jacobi1}.  Because the mass of the $j$th body, which serves as the analog to the circulation $\Gamma_j$, must be positive, condition~\eqref{partial_sum_nonzero} is never an issue in the gravitational problem.

We assign the names $\R_j$ to the final transformed variables and $\kappa_j$ to the transformed circulations. Then $\R_1$ is the displacement from $\r_2$ to $\r_1$ and, similarly, $\R_2$ is the displacement from $\r_3$ to the center of vorticity of the $\r_1$, $\r_2$ subsystem. A similar definition holds for the remaining $\R_j$ with $j<N$, whereas $\R_N$ coincides with our previously defined $\r_0$, the center of vorticity defined in Eq.~\eqref{center_of_vorticity}. Since this quantity is conserved, we have reduced the dimension of the phase space by two.

\subsection{Nambu brackets}
\label{sec:nambu}
The reduced equations of motion we derive will make use of a Nambu bracket, which takes the form 
\begin{equation*}
\poisson{F}{G}_C = -\nabla{C}\cdot\left(\nabla{F} \times \nabla{G} \right),
\end{equation*}
where $C,F,G:\mathbb{R}^3\to\mathbb{R}$, and $C$ is a distinguished function or \emph{Casimir}. The Nambu bracket obeys all the defining properties of a Poisson bracket. Namely, it is a skew-symmetric bilinear operator that obeys the Leibnitz rule and the Jacobi identity. Yoichiro Nambu introduced it in 1973\cite{Nambu:1973} to generalize Hamiltonian mechanics to systems with three-dimensional phase space. Holm et al.'s textbooks provide an excellent overview of the mathematical theory and many applications to problems in optics, classical mechanics, and fluid dynamics\cite{Holm.2011xz9,Holm.2011,Holm.2009}. 

For a system with coordinates $(X,Y,Z)$, the system evolution equations analogous to system~\eqref{dqdt} is
\begin{equation} \label{XYZdot_crossproduct}
\dv{t} \begin{pmatrix} X \\ Y \\ Z \end{pmatrix} = \nabla C \times \nabla H,
\end{equation}
under which any function $F(X,Y,Z)$ evolves according to
\begin{equation*} 
\dv{F(X,Y,Z)}{t} = \poisson{F}{H}_C,
\end{equation*}
in analogy with Eq.~\eqref{Fdot}.

Müller and Névir showed that Gröbli's reduced equations could be reinterpreted as Nambu dynamics\cite{Muller.2014}, but their construction does not solve the problem of those coordinates' singularity. This formulation was subsequently applied to study self-similar collapse in a generalized point-vortex problem\cite{Badin:2018}.

The Nambu formulation of mechanics is often used in situations for which polar coordinates or their Hamiltonian equivalents may be applied. For example Luo et al.\ arrive at an equation equivalent to Eq.~\eqref{generalhamiltonian} below and then introduce polar coordinates. Polar coordinates introduce a singularity at the origin, similar to the singularity of the trilinear coordinate system. The more modern formulation avoids this. This point of view is well articulated in the aptly named lecture notes in Ref.~\onlinecite{Efstathiou:2005}.

\subsection{Application to the three-vortex system}
\label{sec:reduce3vortex}

The Jacobi coordinates and the virtual circulations for the three-vortex problem under assumptions~\eqref{Gamma_signs} and~\eqref{partial_sum_nonzero} are
\begin{equation}\label{jacobi_coordinates}
\begin{aligned}
    \R_1 & =  \r_1 - \r_2; & \kappa_1 & = \frac{\Gamma_1  \Gamma_2}{\Gamma_1 + \Gamma_2}; \\
    \R_2 & = \frac{\Gamma_1\r_1 + \Gamma_2 \r_2}{\Gamma_1 + \Gamma_2} -\r_3; & \kappa_2 & =  \frac{(\Gamma_1 + \Gamma_2)\Gamma_3}{\Gamma_1 + \Gamma_2 + \Gamma_3}; \\
    \R_3 & = \frac{\Gamma_1 \r_1 + \Gamma_2 \r_2 + \Gamma_3 \r_3}{\Gamma_1+\Gamma_2+\Gamma_3}; &  \kappa_3 &= \Gamma_1 + \Gamma_2 +\Gamma_3.
\end{aligned}
\end{equation}
Choosing the center of vorticity $\R_3$ as the origin, we may invert these equations to find
\begin{equation*} 
\begin{aligned}
 \r_1 &= \frac{\Gamma_2}{\Gamma_1 + \Gamma_2}\R_1 +\frac{ \Gamma_3}{\Gamma_1 + \Gamma_2 +\Gamma_3}\R_2; \\
 \r_2 & =\frac{ \Gamma_3}{\Gamma_1 + \Gamma_2 +\Gamma_3} \R_2 - \frac{\Gamma_1}{\Gamma_1 + \Gamma_2} \R_1 ;\\
 \r_3 &= - \qty(\frac{\Gamma_1 + \Gamma_2} {\Gamma_1 + \Gamma_2 +\Gamma_3})\R_2.
\end{aligned}
\end{equation*}

In these coordinates, the Hamiltonian and angular impulses are then
\begin{multline} \label{generalhamiltonian}
  H = -\frac{\Gamma_1 \Gamma_2}{2} \log{\norm{\R_1}^2} -
       \frac{\Gamma_2  \Gamma_3}{2} \log{\norm{\R_2 - \frac{\kappa_1}{\Gamma_2}\R_1}^2} \\
      - \frac{\Gamma_1 \Gamma_3}{2} \log{\norm{\R_2+ \frac{\kappa_1}{\Gamma_1}\R_1}^2},
\end{multline}
and
\begin{equation*} 
  \Theta = \kappa_1 \norm{\R_1}^2 + \kappa_2 \norm{\R_2}^2.
\end{equation*}

For the remainder of the paper, we assume that
\begin{equation*}
\Gamma_1 \ge \Gamma_2 >0.
\end{equation*}
This is generic as two of the vortices must have circulations of matching signs, and we may assume they are positive by reversing the direction of time if necessary. Under this assumption, $\kappa_1>0$, but $\kappa_2$ may take either sign, which plays an essential role in the following analysis.

\subsubsection{The case $\kappa_2>0$}
\label{sec:nambu_positive}

Under assumptions~\eqref{Gamma_signs} and~\eqref{partial_sum_nonzero}, the virtual circulation $\kappa_2$ is positive if $\Gamma_3>0$ or $\Gamma_3<-\Gamma_1-\Gamma_2$. In both these cases, Aref finds that the physical domain $\DomainPhysical$ in the trilinear coordinate system is the interior of an ellipse\cite{Aref.1979}. In the first case, the ellipse lies inside the central triangular region as in Fig.~\ref{fig:trilinear_diagram}(a); in the second, it lies in one of the unbounded regions of the figure.

We normalize the system using Eq.~\eqref{xy_to_pq} and the values of $\kappa_j$ from Eq.~\eqref{jacobi_coordinates}, which gives
\begin{equation*}
Q_1 = \sqrt{\kappa_1} X_1, \,
P_1 = \sqrt{\kappa_1} Y_1, \,
Q_2 = \sqrt{\kappa_2} X_2, \,
P_2 = \sqrt{\kappa_2} Y_2. 
\end{equation*}
Defining $\cR_i = \pair{Q_i}{P_i}$ and $\bcR_i = \pair{Q_i}{-P_i}$, which we can treat as complex variables, the Hamiltonian becomes
\begin{multline*}
H  = -\frac{\Gamma_1  \Gamma_2}{2} \log{\norm{\cR_1}^2}  
- \frac{\Gamma_2 \Gamma_3}{2} \log{\norm{\cR_2 - \frac{\sqrt{\kappa_1 \kappa_2}}{\Gamma_2}\cR_1}^2} \\
- \frac{\Gamma_1  \Gamma_3}{2} \log{\norm{\cR_2 +  \frac{\sqrt{\kappa_1 \kappa_2}}{\Gamma_1}\cR_1}^2} ,
\end{multline*}
and the angular impulse becomes
\begin{equation*}
\Theta  = \norm{\cR_1}^2 + \norm{\cR_2}^2.
\end{equation*}
Both $H$ and $\Theta$ are invariant under the $S^1$ transformation $(\cR_1,\cR_2) \to (e^{i\varphi}\cR_1,e^{i\varphi}\cR_2)$ for arbitrary phase $\varphi$ and depend only on quadratic monomials. Holm suggests the following coordinates for such dynamics\cite{Holm.2011xz9},
\begin{equation}\label{ZXYpositive}
\begin{split}
Z & = \norm{\cR_1}^2 - \norm{\cR_2}^2, \\
X + i Y & = 2 \cR_1 \bcR_2.
\end{split}
\end{equation}
These new coordinates satisfy:
\begin{equation} \label{ThetaZXYPositive}
\Theta^2 = Z^2 + X^2 + Y^2,
\end{equation}
and give a Hamiltonian
\begin{equation*}
\begin{split}
H(X, Y, Z, \Theta) = &-\frac{\Gamma_1  \Gamma_2}{2} \log{\left(\frac{Z+\Theta}{2 \kappa_1} \right)}\\
&-\frac{\Gamma_2  \Gamma_3}{2} \log{\left(\frac{\Theta -Z}{2 \kappa_2} + \frac{\kappa_1(Z+ \Theta) }{2\Gamma_2^2} - \frac{k X}{ \Gamma_2} \right)}\\
&-\frac{\Gamma_1  \Gamma_3}{2} \log{\left(\frac{\Theta -Z}{2 \kappa_2} + \frac{\kappa_1(Z+ \Theta) }{2 \Gamma_1^2} + \frac{k X}{ \Gamma_1} \right)},
\end{split}
\end{equation*}
where $k^2 = \frac{\kappa_1}{\kappa_2}$.

The conservation law~\eqref{ThetaZXYPositive} provides a geometric interpretation of Aref's observation that the physical domain $\DomainPhysical$ is bounded by an ellipse: the natural phase space of the system is the sphere $S^2$. A simple calculation shows that 
\begin{equation*}
Y = - 2 \sqrt{\kappa_1 \kappa_2}\, \R_1 \times \R_2 = (\r_2-\r_1)\times(\r_3-\r_1),
\end{equation*}
so that the great circle $Y=0$, which we will call the equator, corresponds to the set of collinear configurations, i.e., to $\partial\DomainPhysical$ in the trilinear coordinates. In the present coordinate system, the dynamics are regular along this curve. 

It is then an exercise in the chain rule to show that the system evolves under system~\eqref{XYZdot_crossproduct} with $C = 2\Theta^2$ in the coordinates defined by~\eqref{ZXYpositive}. Since $\pdv{H}{Y}=0$, this yields
\begin{equation} \label{XYZdot_positive}
\begin{aligned} 
\dv{X}{t} & = \phantom{-} 4 Y \pdv{H}{Z}; \\
\dv{Y}{t} & = \phantom{-}4 Z \pdv{H}{X}  - 4 X \pdv{H}{Z}; \\
\dv{Z}{t} & = -4 Y \pdv{H}{X}. 
\end{aligned}
\end{equation}

\subsubsection{The case $\kappa_2<0$}
\label{sec:nambu_negative}

The virtual circulation $\kappa_2$ is negative if $-\Gamma_1-\Gamma_2<\Gamma_3<0$. In this, Aref finds that the physical domain $\DomainPhysical$ in the trilinear coordinate system is bounded by a hyperbola\cite{Aref.1979}.

We normalize the system using Eq.~\eqref{xy_to_pq} and the values of $\kappa_j$ from Eq.~\eqref{jacobi_coordinates}, which gives
\begin{equation}\label{XY_to_PQminus}
Q_1 = \sqrt{\kappa_1} X_1, \,
P_1 = \sqrt{\kappa_1} Y_1, \,
Q_2 = -\sqrt{-\kappa_2} X_2, \,
P_2 = -\sqrt{-\kappa_2} Y_2. 
\end{equation}

Defining $\cR_j$ and $\bcR_j$ as in the previous section, we find Hamiltonian
\begin{multline}\label{Hcase2}
H = -\frac{\Gamma_1  \Gamma_2}{2}\log{\norm{\cR_1}^2}  
  -\frac{\Gamma_2  \Gamma_3}{2} \log{\norm{\bcR_2 - \frac{\sqrt{-\kappa_1 \kappa_2}}{\Gamma_2}\cR_1}^2} \\
  -\frac{\Gamma_1  \Gamma_3}{2} \log{\norm{\bcR_2 + \frac{\sqrt{-\kappa_1 \kappa_2}}{\Gamma_1}\cR_1}^2} ,
\end{multline}
and angular impulse
\begin{equation*}
\Theta  = \norm{\cR_1}^2 - \norm{\cR_2}^2.
\end{equation*}

Both $H$ and $\Theta$ are invariant under the $S^1$ transformation $(\cR_1,\cR_2) \to (e^{i\varphi}\cR_1,e^{-i\varphi}\cR_2)$ for arbitrary phase $\varphi$ and depend only on quadratic monomials. Holm suggests the following coordinates for such dynamics\cite{Holm.2011xz9},
\begin{equation}\label{ZXYnegative}
\begin{split}
Z & = \norm{\cR_1}^2 + \norm{\cR_2}^2, \\
X + i Y & = 2 \cR_1 \cR_2.
\end{split}
\end{equation}

These coordinates satisfy
\begin{equation} \label{ThetaZXYminus}
 \Theta^2 = Z^2 - X^2 - Y^2,
\end{equation}
which we know to be conserved. Thus, the trajectory $(X(t), Y(t), Z(t))$ is confined to the upper sheet of a hyperbola of two sheets, which degenerates to a cone when $\Theta=0$. Because $Z\ge0$ by definition, the trajectories lie on the upper sheet. As in the $\kappa_2>0$ case, $Y=0$ when the vortices are collinear, so that Eq.~\eqref{ThetaZXYminus} is a the hyperbola that formed $\partial \DomainPhysical$ in the trilinear coordinates.

The Hamiltonian becomes
\begin{equation}\label{hcase2nambu}
\begin{split}
H(X, Y, Z, \Theta) = & -\frac{\Gamma_1  \Gamma_2}{2}\log{\left(\frac{Z+\Theta}{2 \kappa_1}\right)}   \\
&-\frac{\Gamma_2  \Gamma_3}{2}\log{\left( \frac{Z-\Theta}{2 \kappa_2} + \frac{\kappa_1 (Z+ \Theta)}{2 \Gamma_2^2} - \frac{l X}{\Gamma_2}\right)} \\
&-\frac{\Gamma_1  \Gamma_3}{2}\log{\left( \frac{Z-\Theta}{2 \kappa_2} + \frac{\kappa_1 (Z+ \Theta)}{2\Gamma_1^2} + \frac{l X}{\Gamma_1}\right)},
\end{split}
\end{equation}
where $l^2 = \frac{-\kappa_1}{\kappa_2}$.

The system evolves under equation~\eqref{XYZdot_crossproduct} with $C=2\Theta^2$. Since $\pdv{H}{Y}=0$, this yields 
\begin{equation*} 
\begin{aligned} 
\dv{X}{t} & = \phantom{-}  4 Y \pdv{H}{Z}; \\
\dv{Y}{t} & = -4 Z \pdv{H}{X}  - 4 X \pdv{H}{Z}; \\
\dv{Z}{t} & = -4 Y \pdv{H}{X}.
\end{aligned}
\end{equation*}

\section{The system of three identical vortices}
\label{sec:111}

We illustrate the reduction for the case of three equal circulations $\Gamma_1 = \Gamma_2 = \Gamma_3 = 1$. In this case $\kappa_1 = \frac{1}{2}$ and $\kappa_2=\frac{2}{3}>0$, so the reduction follows Sec.~\ref{sec:nambu_positive} and
\begin{multline*}
H = -\frac{1}{2} \log{(\Theta+Z)}
    -\frac{1}{2} \log{\left(\Theta-\frac{Z}{2}-\frac{\sqrt{3}X}{2}\right)} \\
    -\frac{1}{2} \log{\left(\Theta-\frac{Z}{2}+\frac{\sqrt{3}X}{2}\right)}.
\end{multline*}
This Hamiltonian is unchanged under rotations of the $XZ$ plane by $\pm \frac{2\pi}{3}$, which correspond to permutations of the vortex labels.
The dynamics are equivariant under a rescaling of $\Theta$, so we may take $\Theta=1$ without loss of generality. 

The dynamics are singular at the points on the sphere where the arguments of the logarithms defining $H$ vanish, which occur at three points evenly spaced around the equator, here given by $Y=0$,
$$
(X,Y,Z) = (0,0,-1) \qand (X,Y,Z) = \left(\pm \frac{\sqrt{3}}{2},0,\frac{1}{2} \right).
$$
These are points where two of the three vortices coincide, and the rotational frequencies of the closed orbits surrounding these points diverge as they approach the singular points. Each point on the three meridians running from the north to the south pole through a singularity corresponds to a "tall" isosceles triangle with legs longer than its base.

The system has five equilibria. Three of them lie on the equator,
$$
(X,Y,Z) = (0,0,1) \qand (X,Y,Z) = \left(\pm \frac{\sqrt{3}}{2},0,-\frac{1}{2} \right).
$$
These alternate with the three singular points as one moves around the equator. They correspond to collinear relative equilibria, each with one of the three vortices at the midpoint of a line segment connecting the other two. The other two, which lie at the poles
$$
(X,Y,Z) = (0,\pm 1, 0),
$$
correspond to rigidly rotating equilateral triangular arrangements. Each point on the three meridians running from the north to the south pole through a saddle point on the equator corresponds to a "wide" isosceles triangle whose base is longer than its legs.

Because the dynamics defined by system~\eqref{XYZdot_positive} are regular, the linear stability of all relative equilibria is determined by the eigenvalues of the Jacobian. We include the following elementary calculations to demonstrate their straightforwardness compared to previous formulations of the problem. The Jacobian matrices at, respectively, a triangular relative equilibrium and a collinear one are
$$
J(0,1,0) = \begin{pmatrix} 
  0 & 0 & 3 \\
  0 & 0 & 0 \\
  -3 & 0 & 0
  \end{pmatrix},
\qand
J(1,0,0) = \begin{pmatrix} 
  0 & 3 & 0 \\
  9 & 0 & 0 \\
  0 & 0 & 0
  \end{pmatrix}.
$$
Each has a null eigenvector corresponding to a perturbation in the radial direction, i.e., to a change to the conserved angular impulse $\Theta$. The first has eigenvalues $\pm 3\ii$ and is neutrally stable. The second has eigenvalues $\pm 3\sqrt{3}$ and is a saddle.

\begin{figure}[htbp] 
   \centering
   \includegraphics[width=2in]{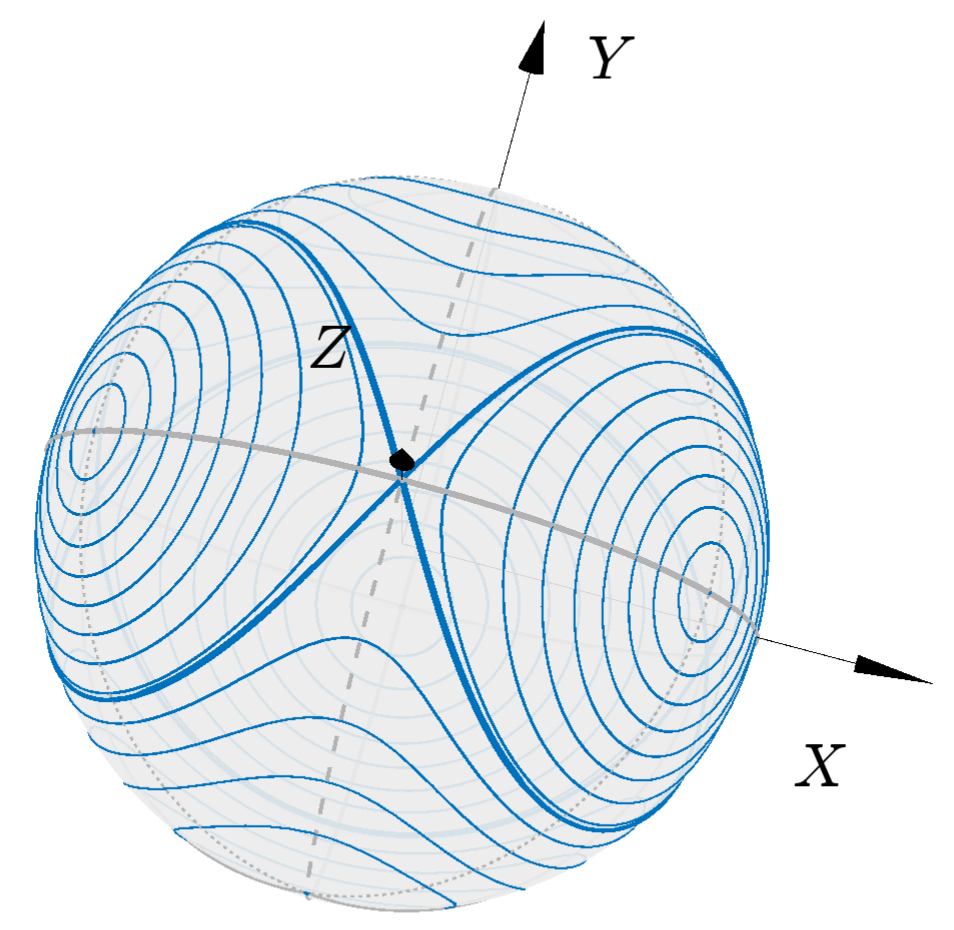} 
   \caption{The phase-sphere of the three-vortex system with identical circulations, plotted with transparency so trajectories on the rear are visible. The equator $Y=0$, where the vortices are collinear, and the meridians, where they form an isosceles triangle, are indicated.}
   \label{fig:111sphere}
\end{figure}

The global phase space for three identical vortices is shown in Fig.~\ref{fig:111sphere}. When this sphere is viewed from above the north pole, it reduces to Aref's phase plane shown in Fig~\ref{fig:trilinear_diagram}(a). The three collinear states are saddle points, and their invariant manifolds coincide in six heteroclinic orbits. The periodic orbits can be classified into two types: two families of orbits that encircle the triangular configurations at the poles and three families that surround the singular points on the equator. 

\begin{figure}[htbp] 
   \centering
   \includegraphics[width=.8 \columnwidth]{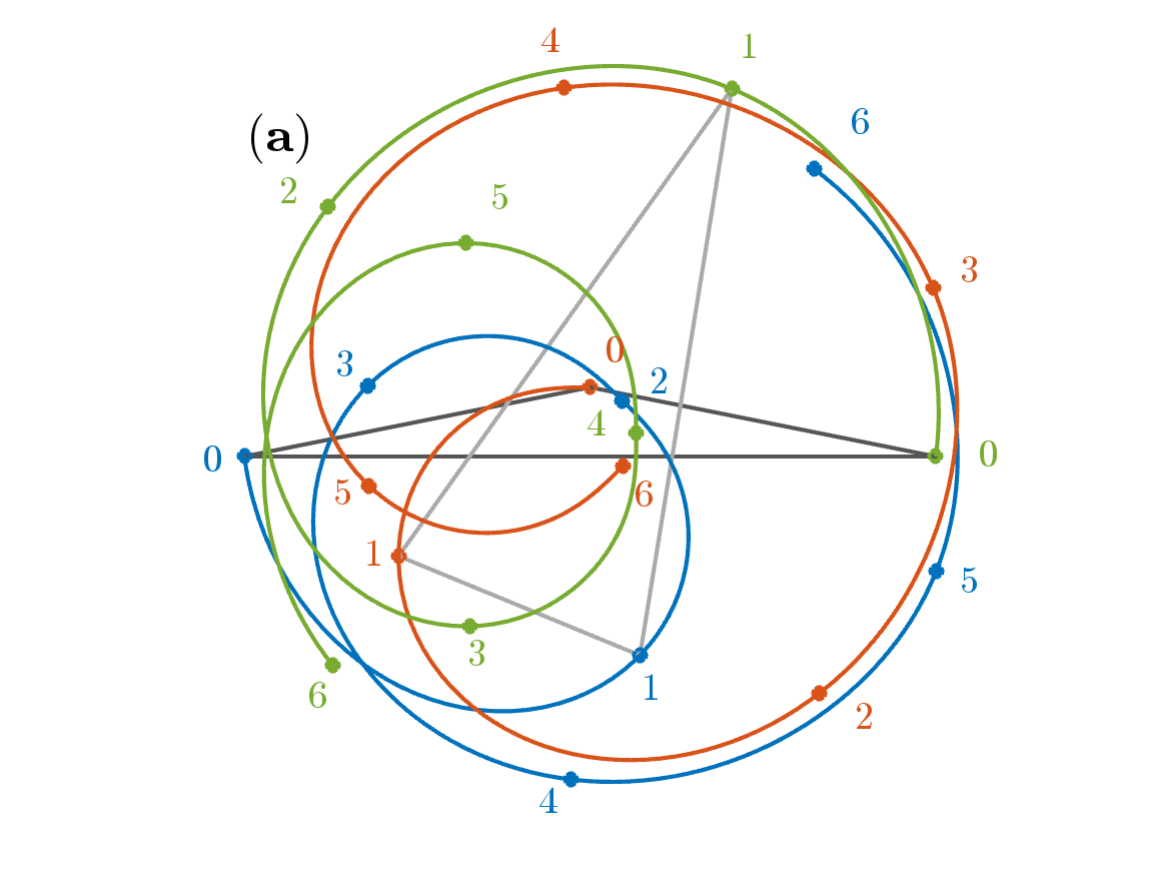} \\
   \includegraphics[width=.8 \columnwidth]{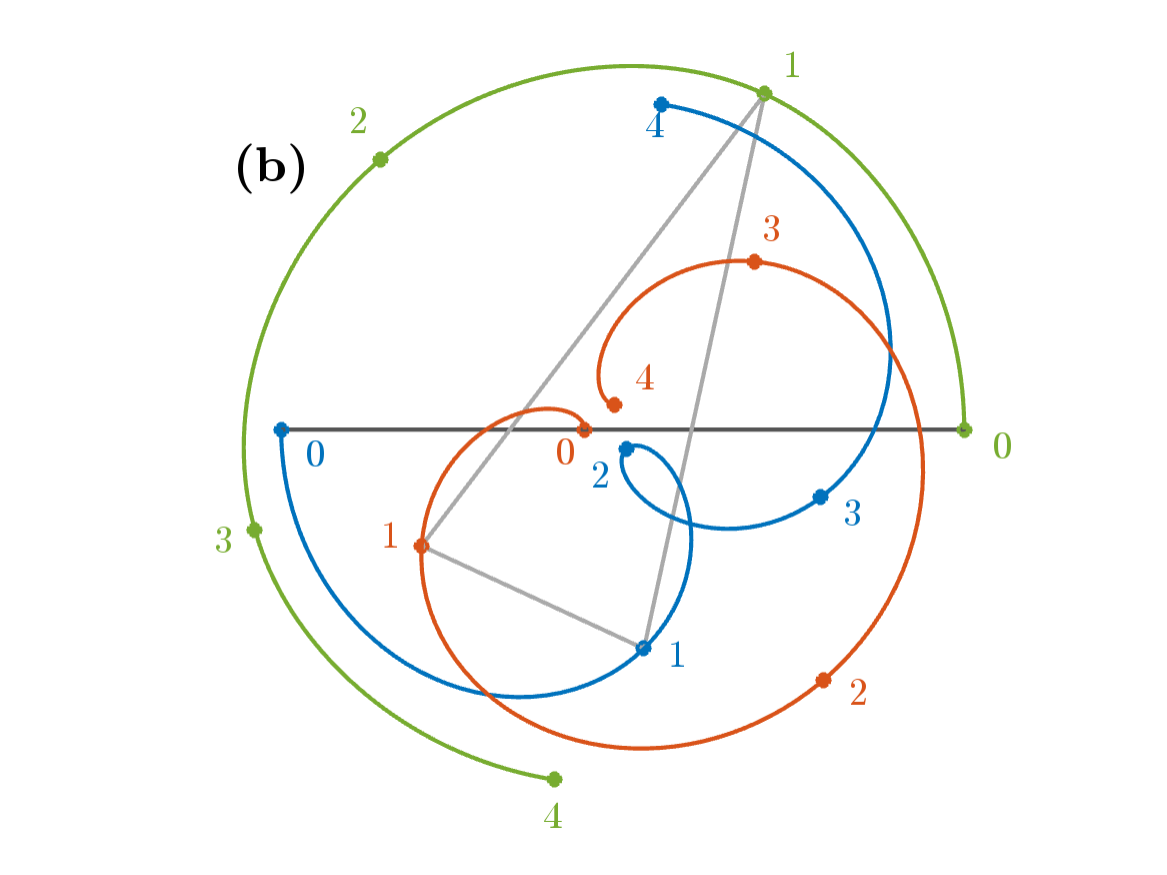} 
   \caption{Periodic orbits near separatrices (relative periodic orbits in physical coordinates). \textbf{(a)} A periodic orbit in the upper hemisphere that makes close approaches to all three saddle points but never crosses the equator. At all the numbered times, corresponding to sixths of a period, the vortices form an isosceles triangle, the first two of which are drawn. \textbf{(b)} A periodic orbit crosses the equator and closely approaches saddle points. At numbered times, corresponding to quarter-periods, the points alternate between collinear and isosceles arrangements.}
   \label{fig:sphere_periodic}
\end{figure}

Fig.~\ref{fig:sphere_periodic}(a) shows a periodic orbit from the family surrounding the north pole with an initial point on the meridian between the $Y$axis and the $Z$ axis in Fig~\ref{fig:111sphere} close to the saddle point. In laboratory coordinates, it is a relative periodic orbit whose initial condition is a "wide" isosceles triangle in which the vortices are nearly collinear. The figure shows one period of motion on the sphere, which crosses all six isosceles meridians but remains in the upper hemisphere.

Fig.~\ref{fig:sphere_periodic}(b) shows a periodic orbit from the family surrounding a singular point on the equator. The initial condition is collinear, with the points nearly equally spaced, corresponding to a point on the equator near a saddle point. Two of the vortices alternate, moving to the center as the orbit approaches two of the saddle points in turn.

\section{Three vortices with circulations $(1,1,-1)$}
\label{sec:11m1}

We let $\Gamma_1 = \Gamma_2 = -\Gamma_3 = 1$, in which case the transformed circulations are
$$
\kappa_1 = \frac{1}{2}, \, \kappa_2 = -2, \, \kappa_3 =1.
$$
Because $\kappa_2<0$, this system is reduced to Nambu form using equations~\eqref{XY_to_PQminus}--\eqref{ZXYnegative}. The Hamiltonian reduces to 
\begin{equation} \label{HXYZ}
  H = -\frac{1}{2} \log{(Z+\Theta)} + \frac{1}{2} \log{\left( Z^2 - X^2\right)}.
\end{equation}
and the angular impulse is given by Eq.~\eqref{ThetaZXYminus}.

The evolution equations are
\begin{subequations}  \label{XYZdot}
\begin{align}
\dv{X}{t} & = \frac{-2Y}{Z+\Theta} + \frac{4 Z Y}{Z^2-X^2}; \label{Xdot} \\
\dv{Y}{t} & = \frac{2X}{Z+\Theta};  \label{Ydot} \\
\dv{Z}{t} & = \frac{4 X Y}{Z^2 - X^2}. \label{Zdot}
\end{align}
\end{subequations}

\begin{figure}[htbp] 
   \includegraphics[width=0.8\columnwidth]{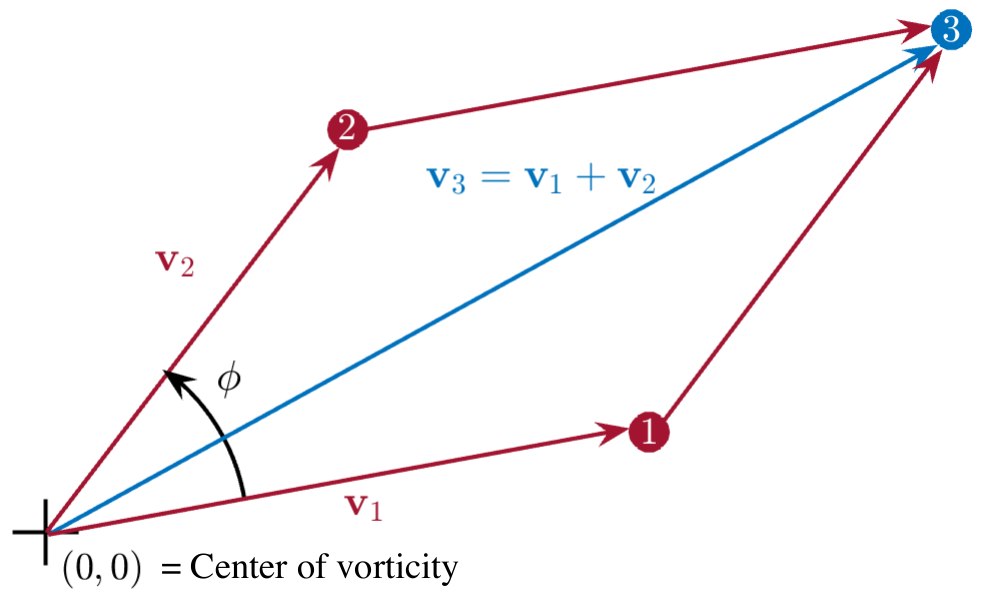} 
   \caption{Diagram used to interpret the $(X,Y,Z,\Theta)$ coordinates. See the text for an explanation. }
   \label{fig:parallelogram}
\end{figure}

It is worth relating the $(X,Y,Z,\Theta)$ coordinate system for this problem to the physical coordinates, and we find a straightforward geometric interpretation. Consider Figure~\ref{fig:parallelogram}. By assumption~\eqref{center_of_vorticity}, the center of vorticity lies at the origin, and we let $\vv_j$ denote the vector from the origin to vortex $j$. Relation~\eqref{center_of_vorticity} implies that $\vv_3 = \vv_1 + \vv_2$ so that the positions of the three vortices and the origin form a parallelogram, a fact mentioned by Gröbli in~[\onlinecite{Grobli:1877},\S3]. We let $\phi$ be the angle from $\vv_1$ to $\vv_2$. We then find by following through the sequence of changes of variables that
\begin{equation} \label{XYZW}
 \begin{split}
       X & = -\norm{\vv_1}^2 + \norm{\vv_2}^2; \\
       Y &= 2 \left(\vv_1 \times \vv_2 \right)\cdot \vb{k} = 2 \norm{\vv_1}\norm{\vv_2} \sin{\phi}; \\
       Z & = \phantom{-}\norm{\vv_1}^2 + \norm{\vv_2}^2; \\
  \Theta & = -2\, \vv_1 \cdot \vv_2 = -2 \norm{\vv_1}\norm{\vv_2} \cos{\phi}.
  \end{split}
\end{equation}

We few observations on these coordinates:
\begin{itemize}
\item $X$ is the signed difference between the lengths of $\vv_1$ and $\vv_2$, so vanishes when the triangle of vortices is isosceles.
\item $Y$ vanishes when the three vortices are collinear, which is not a singularity of the coordinate system.
\item For $\Theta=0$, $\cos{\phi}=0$, so at all times, the three vortices form a right triangle, with vortex 2 at the right angle. Then, trivially, they cannot be collinear, so $Y\neq 0$, which can also be deduced from the singularity of the Hamiltonian~\eqref{HXYZ} when $\Theta=0$.
\end{itemize}

\subsection{Scattering}

\begin{figure*}
\includegraphics[width = 0.8\textwidth]{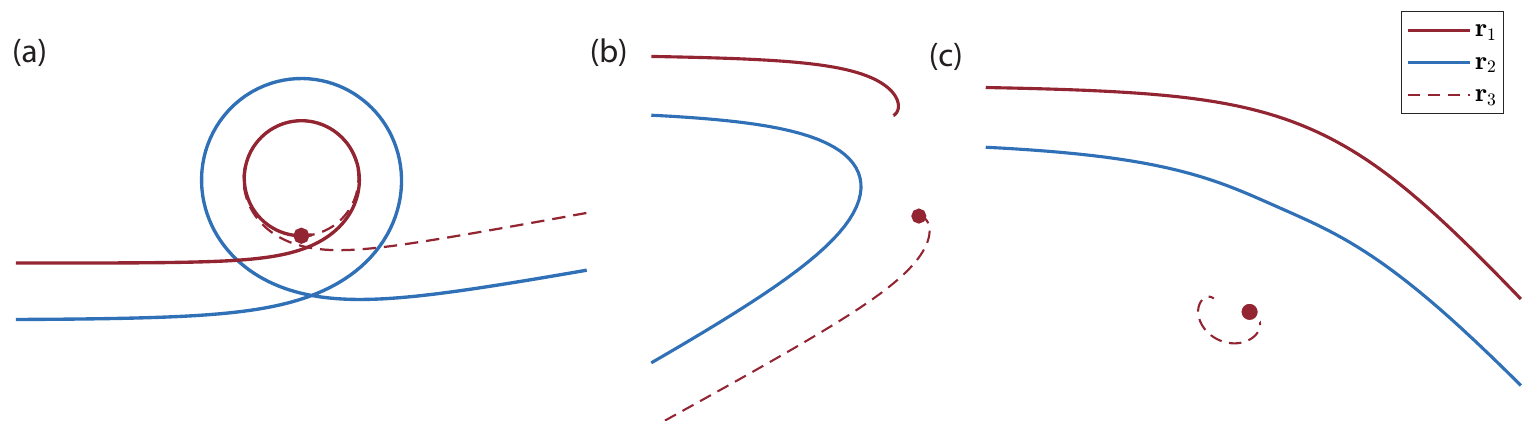}
\caption{\label{fig:scattering_examples}Three solutions of the scattering problem showing (a) Exchange scattering for $\rho=-0.999$. (b) Exchange scattering for $\rho = 2.5$. (c) Direct Scattering for $\rho=3.8$. The vortex dipole arrives from $-\infty$ traveling parallel to the $x$ axis, and vortex 3 sits at rest at the marked point as $t\to-\infty$.}
\end{figure*}

The most noteworthy behavior for this set of circulations is scattering: two vortices with circulations of identical magnitude but opposite orientation form a \emph{dipole} that propagates at constant velocity perpendicular to the line joining them. The presence of a third vortex deflects or \emph{scatters} this motion. Three such scattering solutions are shown in Fig.~\ref{fig:scattering_examples}. While these three solutions obey very similar conditions before the interaction (as the time $t\to-\infty$), their behavior as $t \to \infty$ are quite different. Subfigures (a) and (b) show \emph{exchange scattering} events: the dipole that exits the collision region is not composed of the same two vortices as the dipole that entered. By contrast, subfigure (c) displays \emph{direct scattering}; the same two vortices form the dipole before and after the interaction.

A fundamental question about this scattering is whether a given initial condition leads to direct or exchange scattering. The second question is the change in angle $\Delta\alpha$ between the incoming dipole and the exiting dipole.  Aref derived a formula that answers both questions about scattering, but this is based entirely on integrating the ordinary differential equations and not on the interpretation of the phase diagram in Fig.~\ref{fig:trilinear_diagram}(b)\cite{Aref.1979}. Aref's plot of the dependence of $\Delta\alpha$ on initial contains a sign error that was fixed by Lydon et al.\cite{Lydon.2022}.

\begin{figure}[htbp] 
   \centering
   \includegraphics[width=\columnwidth]{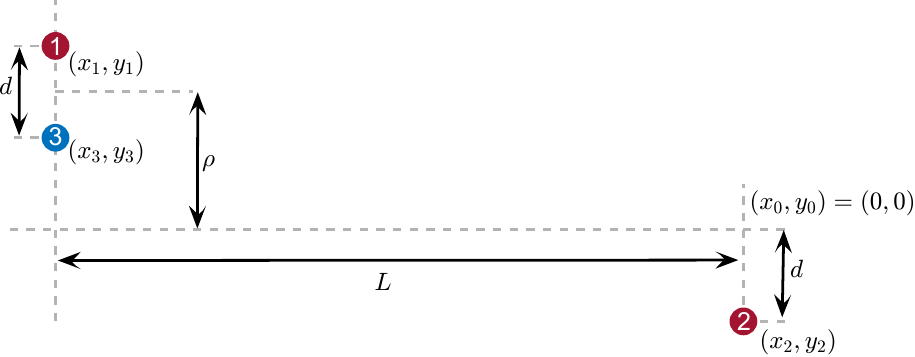} 
   \caption{\label{fig:setup}Setup of the scattering problem. The dipole formed by vortices 1 and 3 propagates toward the target, vortex 2.}
\end{figure}

The reduced system of equations allows us to apply phase space reasoning directly to the scattering problem, so we review it here. A schematic of the scattering experiment is shown in Fig.~\ref{fig:setup}. A dipole consisting of a positive-circulation vortex at position $\r_1 = \pair{-L}{\rho+ \frac{d}{2}}$, where $L\gg1$, and a negative-circulation vortex at position $\r_3 = \pair{-L}{\rho- \frac{d}{2}}$ propagates to the right toward a positive-circulation vortex at position $\r_2 = \pair{0}{-d}$. These are chosen to set $\r_0=0$. Without loss of generality, we take $d=1$. 

Eventually, vortex 3 escapes to infinity as part of a dipole. We call the case when the escaping dipole comprises vortices 1 and 3 a \emph{direct} scattering event and the case when it comprises vortices 2 and 3 an \emph{exchange} scattering event.  

\begin{figure*}[htb] 
   \centering
   \includegraphics[width=0.8\textwidth]{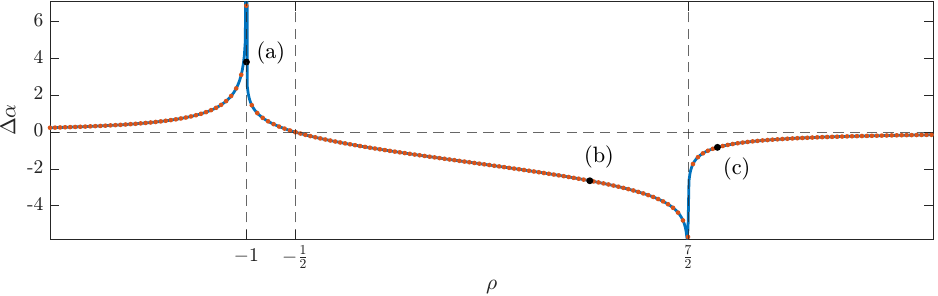} 
   \caption{The deflection of the angle of vortex 3 plotted as a function as the distance $\rho$, showing singularities at $\rho=-1$ and $\rho = \frac{7}{2}$ as expected. The solid line is the result of direct simulation, and the red dots are the formulas derived in Appendix~\ref{sec:appendix}. The points marked (a)--(c) correspond to the three simulations shown in Fig.~\ref{fig:scattering_examples}.}
   \label{fig:theta_vs_rho}
\end{figure*}

Examples are shown in Fig.~\ref{fig:scattering_examples}. The initial conditions are posed as in the schematic, showing exchange scattering in panels (a) and (b) and direct scattering in panel (c). Since vortex 3 has opposite circulation to the two others, it must be part of both the entering and exiting dipoles. We define the \emph{scattering angle} $\Delta \alpha$ as its change of heading; see Eq.~\eqref{alphadot}.  Figure~\ref{fig:theta_vs_rho} shows the scattering angle as a function of the offset $\rho$, with the scattering angles of the three solutions shown in Fig.~\ref{fig:scattering_examples} marked.

If $\abs{\rho} \gg  1$, the isolated vortex will scarcely deflect the dipole, so direct scattering will occur. Previous authors have determined, via fairly involved calculations, that exchange scattering occurs for $-1< \rho< \frac{7}{2}$, and direct scattering outside this interval\cite{Aref.1979, Lydon.2022}. The points $\rho = -1$ and $\rho = \frac{7}{2}$ separate distinct behavior domains in this system, and the scattering angle diverges as $\rho$ approaches these values.

\subsection{Recovering some angles}
\label{sec:action-angle}

Eq.~\eqref{XYZW} shows that the $(X,Y,Z,\Theta)$ variables are insensitive to a rigid rotation of the parallelogram in Fig.~\ref{fig:parallelogram} about the origin and will not allow us to compute the scattering angle. Therefore, we introduce a canonical form of polar coordinates (the action-angle variables of a harmonic oscillator) to recover this angle. 

Returning to the coordinates $\cR_1$ and $\cR_2$ used in Eq.~\eqref{Hcase2}, we let
\begin{equation*}
\begin{aligned}
  \cR_1 = \pair{\sqrt{2I_1} \sin{\theta_1}}{\sqrt{2I_1} \cos{\theta_1}}; \\ 
  \cR_2 = \pair{\sqrt{2I_2} \sin{\theta_2}}{\sqrt{2I_2} \cos{\theta_2}}.
  \end{aligned}
\end{equation*}

Two observations are important here. First, the Hamiltonian depends on the angles only through the combination $\theta_1 + \theta_2$. Second, the vector $\vv_3$ in the figure has argument $\theta_2 - \frac{\pi}{2}$.
Therefore, we  make one additional canonical transformation
\begin{equation*}
\psi_1 = \theta_1 + \theta_2, \quad
\psi_2 = \theta_2, \quad
J_1 = I_1, \quad
J_2 = I_2-I_1.
\end{equation*}
In these variables, the Hamiltonian takes the form (again ignoring additive constants)
\begin{equation*} 
H = \frac{1}{2} \log {\left(4 J_1^2 \sin^2{\psi_1}+4 J_1 J_2 \sin^2{\psi_1}+J_2^2\right)}
-\frac{1}{2} \log {\left( J_1\right)}.
\end{equation*}
Since the equation is cyclic in $\psi_2$, the action $J_2= -\Theta/2$ is conserved. The dynamics of $J_1$ and $\psi_1$ are equivalent to system~\eqref{XYZdot}. We may recover the evolution of $\theta_2=\psi_2$ by integrating
\begin{equation*} 
\dot{\psi}_2 = \frac{2 J_1 \sin^2{\psi_1}+J_2}{4 J_1^2 \sin^2{\psi_1}+4 J_1 J_2 \sin^2{\psi_1}+J_2^2}
\end{equation*}
along a scattering trajectory.
In terms of the Nambu variables, this becomes
\begin{equation*} 
\dot{\theta}_2 = \frac{2 Y^2 \sqrt{\Theta^2+X^2+Y^2}-2 \Theta  X^2}{\left(X^2+Y^2\right) \left(\Theta^2+Y^2\right)}.
\end{equation*}

The angle just calculated describes the argument of $\vv_3$ in Fig.~\ref{fig:parallelogram}, which is distinct from the scattering angle $\alpha=\arg{\dv{z_3}{t}}$ plotted in Fig~\ref{fig:theta_vs_rho}. In terms of the reduced coordinates, we find that
\begin{equation}\label{alphadot}
\dv{\alpha}{t} = -\frac{8 \Theta  Y^2}{\left(X^2+Y^2\right) \left(\Theta^2+Y^2\right)}.
\end{equation}
Integrating this over a trajectory then gives $\Delta\alpha$. This calculation is described in Appendix~\ref{sec:appendix}. It is equivalent to a calculation by Lydon and is included for completeness\cite{Lydon.2022}.

\subsection{Phase space of the $(1,1,-1)$ system}
\label{sec:phase_space_1_1_m1}

We first derive the fixed points and singularities of system~\eqref{XYZdot} before visualizing the system's phase space. We set the right-hand sides of system~\eqref{XYZdot} to zero while enforcing the constraints~\eqref{ThetaZXYminus} and  $Z \ge 0$. Similarly, we find singularities where the argument of either logarithmic term in the Hamiltonian~\eqref{HXYZ} vanishes, enforcing the same two constraints. Which equilibria and singularities exist depends on $\Theta$. 

When $\Theta<0$, the system has two equilibria $\EqTriPM$ and a singularity $\SingPP$ found by setting  $Z+\Theta=0$, which requires $X=Y=0$.  These are
\begin{equation*} 
    \EqTriPM  = 
    \begin{pmatrix} 0 \\ \pm\sqrt{3}\Theta \\ -2\Theta \end{pmatrix} \qand
    \SingPP = \begin{pmatrix} 0 \\ 0 \\ -\Theta \end{pmatrix}.
\end{equation*}

When $\Theta=0$, there are no equilibria, but the system is singular when $Z=\abs{X}$, which requires $Y=0$.

When $\Theta>0$, the system has a single equilibrium
\begin{equation*} 
\EqMinus \equiv \begin{pmatrix} X_0 \\ Y_0 \\ Z_0 \end{pmatrix} = \begin{pmatrix} 0 \\ 0 \\ \Theta \end{pmatrix}
\end{equation*}
and no singularities.

The fixed points $\EqTriPM$ and $\EqMinus$  are \emph{relative equilibria} in the laboratory coordinates, i.e., they are equilibria when viewed in an appropriate rotating reference frame. We may interpret them using Eq.~\eqref{XYZW}. For both equilibria $X=0$ implies $\norm{\vv_1}=\norm{\vv_2}$. The equilibrium $\EqTriPM$ exists for $\Theta<0$. For the equilibrium $\EqTriPM$, the value of the component $Z = -2 \Theta$ implies that $\phi = \pm \frac{\pi}{3}$ and the three vortices lie at the vertices of an equilateral triangle, motivating the naming convention.  The equilibrium $\EqMinus$ exists for $\Theta>0$. This implies $\phi = \pi$ so that the three vortices are collinear with the two positive vortices equally spaced from the negative vortex at the center.  The subscript $-1$ indicates that the vortex with circulation $-1$ sits at the center. By similar reasoning, we find that at the singularity $\SingPP$, the two vortices with circulation $+1$ coincide, again motivating the notation.

\begin{figure*}[htbp] 
   \centering
   \includegraphics[width=0.8\textwidth]{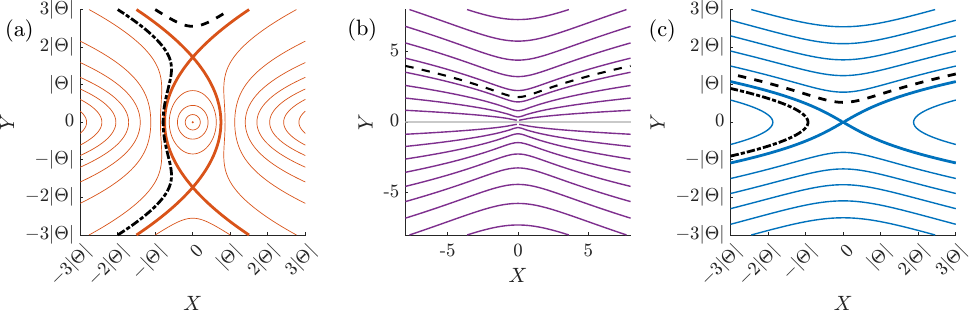}
   \caption{\label{fig:threephaseplanes}The $XY$ phase planes of system~\eqref{XYZdot}. (a) The case $\Theta<0$ with singularity $\SingPP$ (point) and triangular configurations at the intersections of the thick curves. (b) The case $\Theta=0$. The gray line $Y=0$ is singular. (c) the case $\Theta > 0$ with collinear equilibrium $\EqMinus$.  Note that the contours are not evenly-spaced level sets of the energy~\eqref{HXYZ} but were chosen to illustrate the topology clearly.} 
\end{figure*}

By the conservation law~\eqref{ThetaZXYminus}, the phase space of system~\eqref{XYZdot} is the upper sheet of a two-sheeted hyperbola. We visualize the dynamics by projecting this surface into the $XY$ plane in Fig.~\ref{fig:threephaseplanes}. The conserved angular impulse $\Theta$ is a bifurcation parameter, but up to scaling when $\Theta\ne0$, there are only three possible phase planes.

For $\Theta<0$ in panel (a), the point at the origin is the singularity $\SingPP$. The two equilibria $\EqTriPM$ sit on the $Y$ axis and are saddle points connected by a pair of homoclinic orbits. The two homoclines surround a family of periodic orbits, which shrink to a point at $\SingPP$. Each corresponds to a hierarchical orbit in which the two positive vortices orbit about each other rapidly while their mutual center of vorticity and the third vortex orbit each other; Gröbli computed this orbit in closed form and plotted it in Ref.~[\onlinecite{Grobli:1877}, Fig. 1]. As the diameter of these closed orbits goes to zero, the rotation rate of this tightly bound pair diverges, and the orbits approach the singularity $\SingPP$. The unbounded portions of the stable and unstable manifolds separate the remainder of the phase plane into four unbounded quadrants. This will be important for the scattering problem.

When $\Theta=0$ in panel~(b), the entire $X$ axis is singular, and all solutions are confined to the upper or lower half-planes. 
For $\Theta>0$ in panel~(c), the collinear equilibrium $\EqMinus$ at the origin is a saddle point. Its invariant manifolds also separate the plane into four unbounded quadrants. 

The phase plane for $\Theta<0$ in panel~(a) corresponds to the upper disconnected component of $\DomainPhysical$ in Fig.~\ref{fig:trilinear_diagram}(b), the phase plane for $\Theta=0$ in panel~(b)  to Fig.~\ref{fig:trilinear_diagram}(c), and  phase plane for $\Theta>0$ in panel~(c) to the lower disconnected component of $\DomainPhysical$ in Fig.~\ref{fig:trilinear_diagram}(b).

Finally, note that rescaling $X$ and $Y$ by $\abs{\Theta}$ and $t$ by $\abs{\Theta}^{-1}$ (when $\Theta\ne0$) shows that the dynamics for any negative (respectively, positive) value of $\Theta$ has a phase plane equivalent to that shown in panel (a) (respectively, panel (c)).

\subsection{Explaining the scattering}
\label{sec:scattering}

The analysis of the previous section enumerated the ingredients needed to explain the behavior of the three-vortex scattering problem set up in Fig.~\ref{fig:setup}. The most important features of a phase plane in organizing the dynamics are the invariant sets: equilibria, periodic orbits, and their stable and unstable manifolds which form separatrices between regions of this behavior. 
The goal of this section is to show that the transitions between direct and exchange scattering at $\rho = -1$ and $\rho = \frac{7}{2}$ in Fig.~\ref{fig:theta_vs_rho} are due to these features.

The separatrices shown in Fig.~\ref{fig:threephaseplanes} divide the phase plane into families of trajectories with identical topology, and the topology of the phase plane is determined, in turn, by the conserved parameter $\Theta$. Panel (a) depicts the case $\Theta<0$, where the energy level of the separatrices equals that of the rotating triangular configurations $\EqTriPM$, given by
\begin{equation}\label{E_Xpm}
E(\EqTriPM)= \frac{1}{2} \log{(-4\Theta)}.
\end{equation}

The energy in the two regions to the left and right of $\EqTriPM$ (those containing the $X$-axis) is lower than $E(\EqTriPM)$, while the energy in the regions above and below the separatrices is higher than $E(\EqTriPM)$.

Panel (c) shows the case $\Theta>0$, where the energy level on the separatrices equals that of the collinear equilibrium $\EqMinus$, which we compute to be 
\begin{equation}\label{E_X0}
E(\EqMinus)= \frac{1}{2} \log{\frac{\Theta}{2}}.
\end{equation}
The energy in the two regions to the left and right of $\EqMinus$ (those containing the $X$-axis) is lower than $E(\EqMinus)$, while the energy in the regions above and below the separatrices is higher than $E(\EqMinus)$.

We must compare these energies with those of the pre-scattering condition depicted in Fig.~\ref{fig:setup}. In this arrangement, the center of vorticity is at the origin, so we may compute the limiting behavior of $X$ and $Y$ using the equations in~\eqref{XYZW}.  We directly compute that, independent of $L$, 
\begin{equation} \label{Theta_rho}
\Theta = 1 + 2 \rho.
\end{equation}
We assume that as $t \to -\infty$, $L \to \infty$, thus $\norm{\vv_1} \to \infty$ while $\norm{\vv_2}$ is finite, so that $X \to -\infty$, this also implies that $Y \to \infty$. Thus, for the situation depicted in Fig.~\ref{fig:setup}, trajectories in the phase planes depicted in Fig.~\ref{fig:threephaseplanes} arrive from infinity from the northwest direction heading southeast. 

Then, suppose the initial energy exceeds the separatrix energy. In that case, the trajectory begins above the separatrix and crosses the line $X=0$, where $\dv{Y}{t}=0$ before escaping to infinity in the northeast direction. Because $X\to+\infty$ as $t\to\infty$, $\norm{\vv_2}$ must diverge, and this is an exchange scattering event. At the instant the solution crosses $X=0$, then $\norm{\vv_1}=\norm{\vv_2}$ at which point the vectors $\vv_1$ and $\vv_2$ form the legs of an isosceles triangle.

If the initial energy lies below the separatrix energy, then the trajectory begins below the separatrix. It will cross $Y=0$ at which point $\dv{X}{t}$ =0. When $Y=0$, $\sin{\phi}=0$, and the three vortices are collinear. Because $X<0$ along the entire trajectory and $X\to-\infty$ as $t\to\infty$, then $\norm{\vv_1}\to \infty$ and this solution represents a direct scattering event.

For $\Theta<0$, that is, for $\rho < - \frac{1}{2}$, the critical energy is given by Eq.~\eqref{E_Xpm}, which, combined with Eq.~\eqref{Theta_rho} gives a critical energy
\begin{equation} \label{rho_critical_minus}
\rhocritminus = -1.
\end{equation}

For $\Theta>0$, that is, for $\rho > - \frac{1}{2}$, the critical energy is given by Eq.~\eqref{E_X0}, which, combined with Eq.~\eqref{Theta_rho} gives a critical energy
\begin{equation} \label{rho_critical_plus}
\rhocritplus = \frac{7}{2}.
\end{equation}
Fig.~\ref{fig:theta_vs_rho} shows the deflection in the angle of vortex two following the interaction is singular as $\rho\to-1$ and $\rho \to \frac{7}{2}$ as expected.

To calculate the scattering angle, we must integrate Eq.~\eqref{alphadot} over each scattering trajectory. This is equivalent to a calculation by Lydon et al.\cite{Lydon.2022}, and we present it for completeness in Appendix~\ref{sec:appendix}.

We end this section by remarking that the values $\rho = 1$ ($\Theta=-1$), $\rho=-\frac{1}{2}$ ($\Theta=0$), and $\rho=\frac{7}{2}$ ($\Theta=8$) divide the space of initial conditions into four intervals on which the behavior is qualitatively distinct. Gröbli made the same observation (using a constant $\lambda=\Theta/2$) as did Lydon et al.\cite{Grobli:1877,Lydon.2022}, but without referencing a phase plane to organize the orbits. Because both prior works focus on integrating the ODE system via quadrature, these intervals are distinguished mainly by the change in the algebraic forms of those integrals rather than the phase space topology.

\subsection{The borderline case $\Theta=0$}

The approach taken here is especially illuminating for the transition at $\Theta=0$ where Lydon noticed an algebraic change in the form of the integrals but found no visible discontinuity in the scattering angle in Fig.~\ref{fig:theta_vs_rho}\cite{Lydon.2022}.

As $\rho$ increases from $-\infty$ to $\infty$, it crosses the two critical values found above and, in between them, crosses $\rho = -\frac{1}{2}$ at which point $\Theta=0$. In this case, the conservation law~\eqref{XYZW} confines the dynamics to a cone, whose projection into the $XY$ plane is shown in Fig.~\ref{fig:threephaseplanes}(b)

The critical values of $\rho$ separating exchange scattering from direct scattering are those for which the energy of the initial condition as $L \to \infty$ equals that of the separatrices, i.e., the energy of the hyperbolic fixed points. The singular case $\Theta=0$, when the conservation law~\eqref{XYZW} confines the dynamics to a cone, therefore corresponds to $\rho = -\frac{1}{2}$. 

\begin{figure}[htbp] 
   \centering
   \includegraphics[width=.45\textwidth]{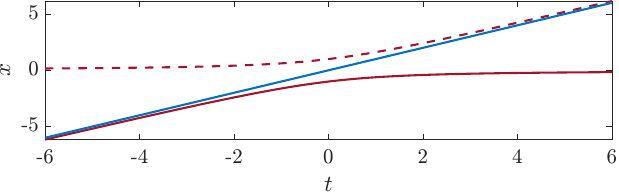}
   \caption{\label{fig:rho_minus_half}The $x$-component of the solution for $\rho=0$, corresponding to a trajectory in the middle phase plane of Fig.~\ref{fig:threephaseplanes}. }
   \end{figure}

The schematic in Fig.~\ref{fig:setup}, which is defined for finite $L$, is somewhat misleading, as the trajectories of all three vortices lie along straight lines parallel to the line connecting vortices 1 and 3 in the figure and are not horizontal. Rotating the coordinate system so that the trajectories are horizontal, we find that 
\begin{align*}
x_1 &= \frac{t-\sqrt{t^2+4}}{2}, & x_2 &= \frac{t+\sqrt{t^2+4}}{2}, & x_3 & = t;\\
y_1 &= -1,& y_2 &= -1,& y_3 &= -2.
\end{align*}
The dynamics of this case are shown in Fig.~\ref{fig:rho_minus_half} and were known to Gröbli\cite[\S4]{Grobli:1877}. Vortex 1 slows down and comes to rest at $x=0$, transferring its energy to vortex 2.

\section{Generalization to $\Gamma_2 \neq 1$}
\label{sec:generalization}

In this section, we will generalize to the case in which $\Gamma_1 = -\Gamma_3 = 1$ but $0<\Gamma_2 = \Gamma \neq 1$, i.e., the case when the remaining vortex has a distinct positive circulation. Exchange scattering is no longer possible since vortices 2 and 3 can no longer form a dipole and escape. Consequently, some bifurcation must reconfigure the phase plane dynamics depicted in Fig.~\ref{fig:threephaseplanes}. 

The Jacobi coordinate reduction yields transformed circulations 
$$
\kappa_1 =\frac{\Gamma}{\Gamma + 1}, \, 
\kappa_2  = -\frac{1+\Gamma}{\Gamma}, \, 
\kappa_3 = \Gamma.
$$
As in the previous case $\kappa_2<0$, so change of variables and Hamiltonian structure of Sec.~\ref{sec:nambu_negative} apply, and the Hamiltonian~\eqref{hcase2nambu} reads
\begin{equation*}
\begin{split}
H(X,Y,Z,\Theta)
=& \frac{\Gamma}{2} \log({Z(\Gamma^2+1)+(1-\Gamma^2)\Theta-2\Gamma X}) \\
 &- \frac{\Gamma}{2} \log({{Z+ \Theta}})
 + \frac{1}{2} \log(Z+X).
\end{split}
\end{equation*}

While we were unable to find as useful a geometric interpretation of the coordinates as in Eq.~\eqref{XYZW}, we still have that $Y=0$ whenever the three vortices are collinear. Moreover $X\to-\infty$ as $\norm{\r_2-\r_3} \to \infty$ and $X\to+\infty$ as $\norm{\r_1-\r_2} \to \infty$. This last observation allows us to discriminate between direct and exchange scattering. 

The system evolves according to
\begin{equation} \label{XYZGdot}
\begin{split} 
\dv{X}{t} & = \frac{2 \Gamma Y}{Z+ \Theta}  
  -  \frac{2 Y}{X+Z} -
  \frac{2  \Gamma  ( 1+ \Gamma^2)Y}{(\Gamma^2+1)Z+(1-\Gamma^2)\Theta-2\Gamma X}; \\
\dv{Y}{t} & = 2 - \frac{2  \Gamma X}{Z+ \Theta} + 
  \frac{2 \Gamma(1 + \Gamma^2)X - 4\Gamma^2 Z }{Z(\Gamma^2+1)+(1-\Gamma^2)\Theta-2\Gamma X} ; \\
\dv{Z}{t} & = \frac{2 Y} {Z+X}- \frac{4 \Gamma^2 Y}{ Z(\Gamma^2+1)+(1-\Gamma^2)\Theta-2\Gamma X}. 
\end{split}
\end{equation}
 
\subsection{The phase space for $\Gamma_2 \neq 1$}
\label{sec:phase_space_generalized}
The equilibria of~\eqref{XYZGdot}, which must satisfy both $Z >0$ by~\eqref{ZXYnegative} and satisfy the constraint~\eqref{ThetaZXYminus} are:

\begin{equation} \label{Gequlilbria}
\begin{split}
\EqTriPM &=  \Theta
  \begin{pmatrix} 
    \frac{\Gamma (\Gamma-1)}{\Gamma+1}  \\ 
    \pm \sqrt{3} \Gamma 
  \end{pmatrix}, 
  \Theta < 0 ;\\ 
\EqMinus &=  \frac{\Gamma \Theta}{\Gamma^2 -1}
  \begin{pmatrix}
      1-2\Gamma^2 + \sqrt{4 \Gamma^2 -3} \\ 
     0 
  \end{pmatrix}, \Theta >0, 
   \Gamma >\frac{\sqrt{3}}{2}; \\ 
\EqGamma &= \frac{\Gamma \Theta}{\Gamma^2 -1}
  \begin{pmatrix} 
    1-2\Gamma^2 - \sqrt{4 \Gamma^2 -3} \\ 
    0 
  \end{pmatrix}, 
   \Theta >0,  \frac{\sqrt{3}}{2} < \Gamma < 1 ; \\ 
\EqPlus &= \frac{\Gamma \Theta}{\Gamma^2 -1}
  \begin{pmatrix} 
  1-2\Gamma^2 - \sqrt{4 \Gamma^2 -3} \\
   0 
  \end{pmatrix}, 
    \Theta <0,  \Gamma > 1. 
    \end{split}
\end{equation}

The same formula describes these last two, but they represent different vortex configurations and are defined for different parameter values. Only the $X$ and $Y$ coordinates are displayed; the $Z$ coordinate is the positive solution to~\eqref{ThetaZXYminus}.

The system also has singularities at
\begin{equation} \label{Gsingularities}
\begin{split} 
\SingPGamma &= 
  \begin{pmatrix} 0  \\ 0 
  \end{pmatrix}, \Theta < 0, \Gamma >0 ; \\
\SingMGamma &= \begin{pmatrix} \frac{2 \Gamma \Theta}{\Gamma^2 -1}  \\ 0 
\end{pmatrix},  \Theta(\Gamma-1)>0.
\end{split}
\end{equation}

\begin{figure}[htbp] 
   \centering
\includegraphics[width=0.45\textwidth]{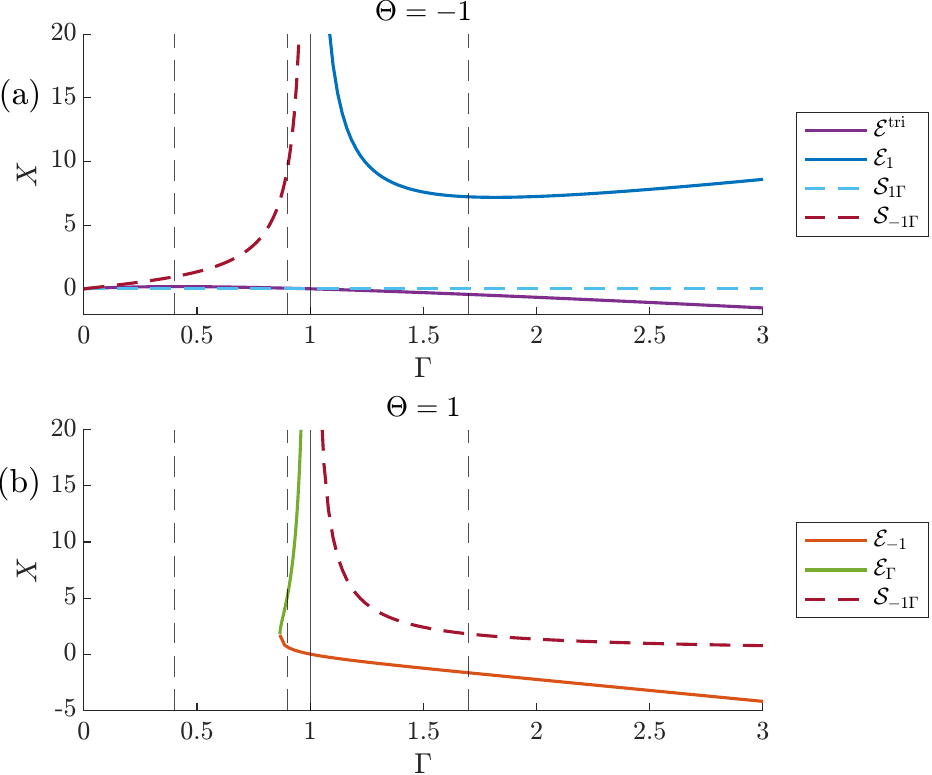} 
\caption{\label{fig:bifurcations_in_Gamma}The $X$ component of the equilibria (solid lines) and singularities (dashed lines) given in Eqs.~\eqref{Gequlilbria} and~\eqref{Gsingularities} for \textbf{(a)}~$\Theta=-1$ and \textbf{(b)}~$\Theta=1$. Figs.~\ref{fig:threephaseplanes}, \ref{fig:gammabigone}, \ref{fig:gammainsqrt1}, and~\ref{fig:gammaless0.8} show phase plane diagrams at the $\Gamma$ values indicated by the vertical lines.}
   
\end{figure}

The dependence of the equilibria and singularities on $\Gamma$ and $\Theta$ are most easily understood graphically using a bifurcation diagram, as shown in Fig.~\ref{fig:bifurcations_in_Gamma}. Only the equilibria $\EqTriPM$ and $\EqMinus$ and the singularity $\SingPGamma$ exist for $\Gamma=1$ and satisfy $\SingPGamma \to \SingPP$ as $\Gamma\to 1$. The other equilibria and singularities all satisfy $Y=0$ and diverge with $X\to +\infty$ as $\Gamma\to\pm 1$. The points $\SingMGamma$, $\EqPlus$, and $\EqGamma$ all diverge to $\infty$ as $\Gamma\to 1^{\pm}$. The equilibria $\EqMinus$ and $\EqGamma$ merge in a saddle-node bifurcation at $\Gamma = \frac{\sqrt{3}}{{2}}\approx 0.866$.

The equilibria $\EqGamma$ and $\EqPlus$ correspond to collinear arrangements with the vortices of strength $\Gamma$ and $1$ in the middle, respectively.  The singularity $\SingMGamma$ corresponds to the limit of a family of hierarchical orbits in which vortices $2$ and $3$, with circulations $\Gamma$ and $-1$, form a tight pair orbiting vortex $1$ some distance away. 

\begin{figure}[htb]
    \centering
     \includegraphics[width= 0.45\textwidth]{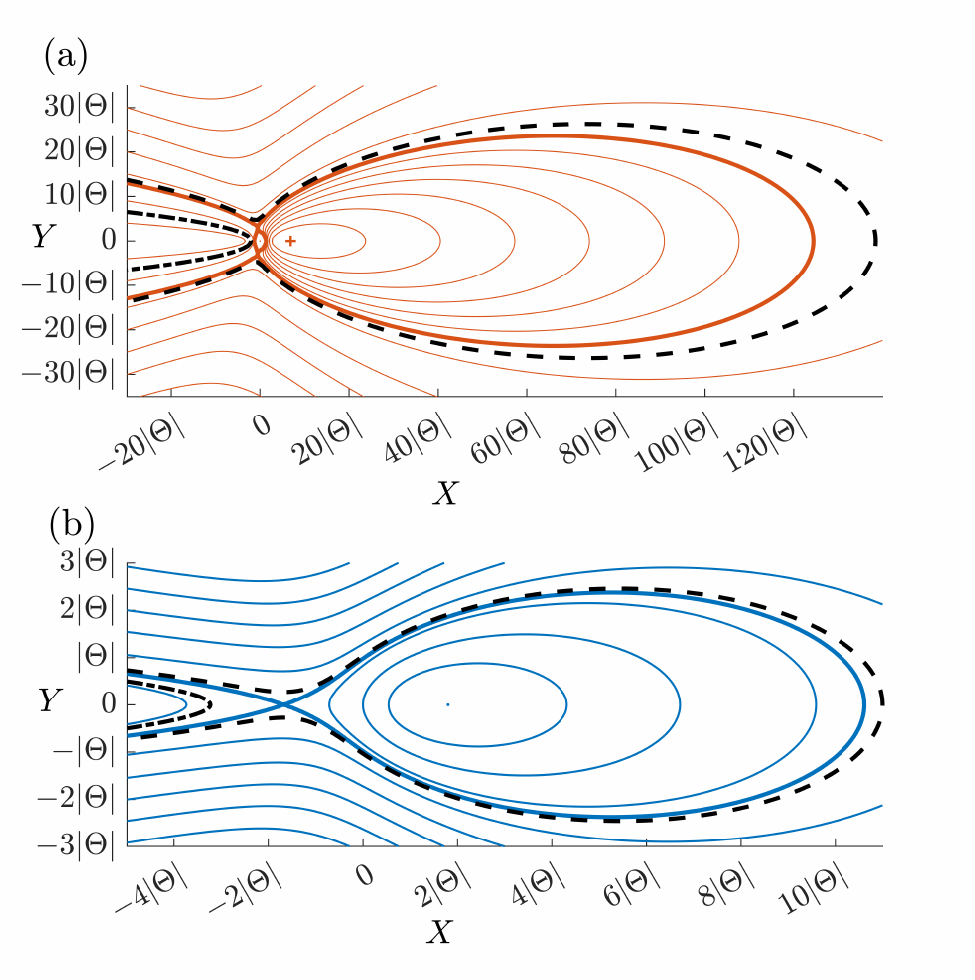}
    \caption{\label{fig:gammabigone}The phase planes for $\Gamma = 1.7 > 1$.  \textbf{(a)} $\Theta <0$,  showing the equilibria $\EqTriPM$ at the separatrix intersections, the singular point $\SingPGamma$ (point) and the collinear state $\EqPlus$ (+).
    \textbf{(b)}~$\Theta > 0$, with singular point $\SingMGamma$ (point) and collinear state $\EqMinus$ at the separatrix intersection. }
    
\end{figure}

We now consider the phase space as $\Gamma$ varies, again plotting the projection of the upper sheet of the hyperboloid onto the $XY$ plane. First, we show the case $\Gamma>1$ as shown in Fig.~\ref{fig:gammabigone}. For $\Theta<0$,  a collinear state $\EqPlus$ appears on the $X$-axis to the right of the region of closed orbits seen in Fig.~\ref{fig:threephaseplanes}, while for $\Theta>0$ a new singular state $\SingMGamma$ appears on the positive $X$-axis. Each of these is surrounded by a family of periodic orbits that limit to a separatrix. In contrast with Fig.~\ref{fig:threephaseplanes}, all orbits in the right half plane cross the $X$-axis and do not extend to $\infty$. 

\begin{figure}[htb]
    \centering
     \includegraphics[width= 0.45\textwidth]{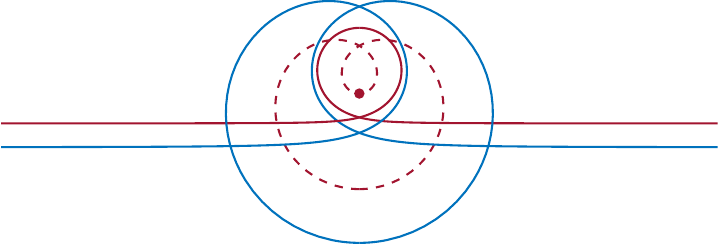}
    \caption{An extended direct scattering solution with $\Gamma=2$ and $\rho=-\frac{9}{2}$. This is a direct simulation of a solution whose trajectory Gröbli computed in closed form.}
    \label{fig:grobliFig5}
\end{figure}

The family of unbounded orbits corresponding to \emph{exchange scattering} has been replaced by a family of orbits that cross the $X$-axis and approach infinity heading southwest. We call these \emph{extended direct scattering} orbits. One such orbit with $\Gamma=2$ is shown in Fig.~\ref{fig:grobliFig5}. Remarkably, The coordinates of this trajectory, but not its time-parameterization, are given by Gröbli and displayed in his dissertation [\onlinecite{Goodman:2024}, Eqns (7.17), (7.19), and (7.20), and Fig. 5]. As $t \to \pm \infty$, vortices 1 and 3 form a dipole that moves along a nearly straight, while at intermediate times, vortex 3 has changed partners and forms a dipole with vortex 2 that along a roughly circular orbit.

\begin{figure}[htb]
    \centering
     \includegraphics[width= 0.47\textwidth]{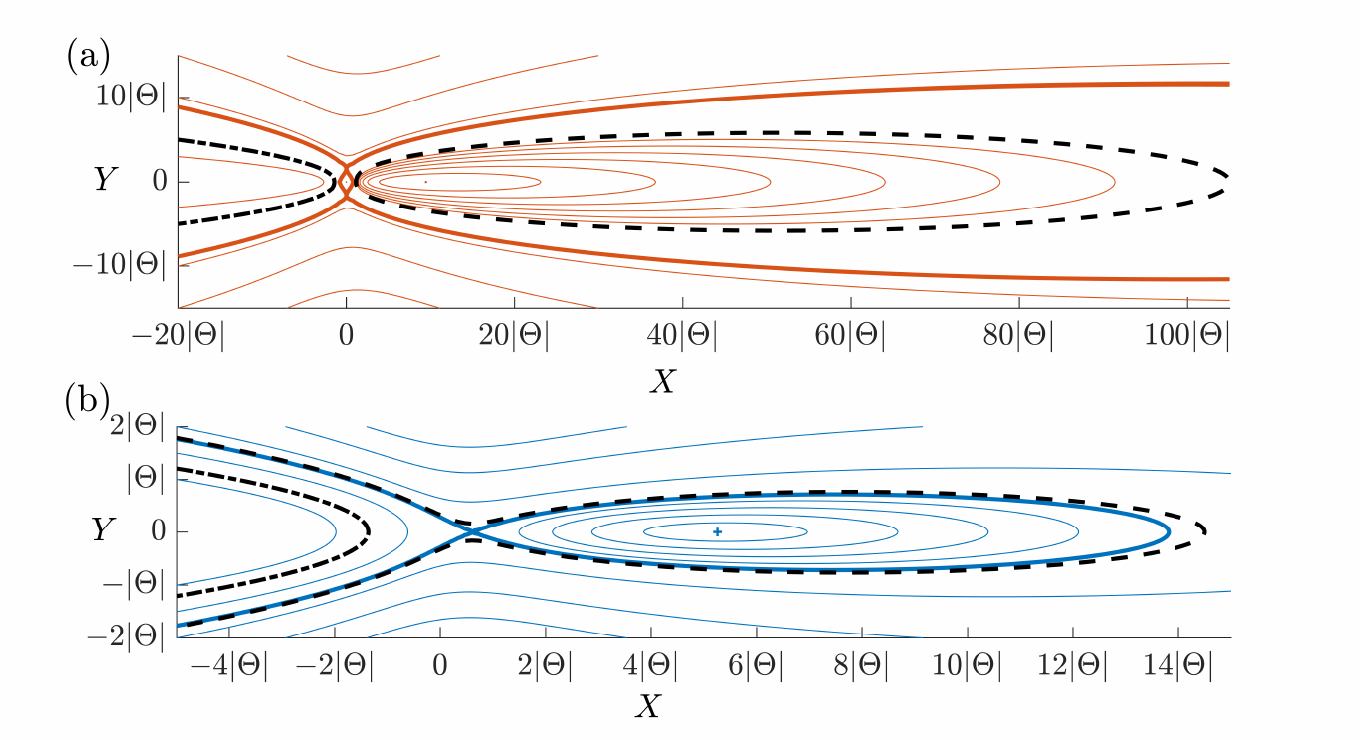}
    \caption{The phase planes for $\Gamma = 0.9\in \qty( \frac{\sqrt{3}}{2}, 1)$. \textbf{(a)}  $\Theta < 0$, with two singular points, $\SingPGamma$ and $\SingMGamma$ (points), and two equilibria $\EqTriPM$ at the separatrix intersections. \textbf{(b)}   $\Theta >0$, with collinear equilibria $\EqGamma$ (+) and $\EqMinus$ at the separatrix intersection.}
    \label{fig:gammainsqrt1}
\end{figure}

Fig.~\ref{fig:gammainsqrt1} shows representative phase planes with $\frac{\sqrt{3}}{2}<\Gamma<1$. While the topology in Figs.~\ref{fig:gammabigone} and~\ref{fig:gammainsqrt1} looks the same, they differ in the kinds of singularities and fixed points.
For $\Theta < 0,$  the singularity $\SingPGamma$ remains unchanged from Fig.~\ref{fig:gammabigone}, while the equilibrium $\EqPlus$ to the right of the origin is replaced by a singularity $\SingMGamma$. 
For $\Theta >0$, the equilibrium $\EqMinus$ is unchanged from Fig.~\ref{fig:gammabigone}, while the singular point $\SingMGamma$ is replaced by the equilibrium $\EqGamma$.  

\begin{figure}[htbp]
    \centering
     \includegraphics[width=\columnwidth]{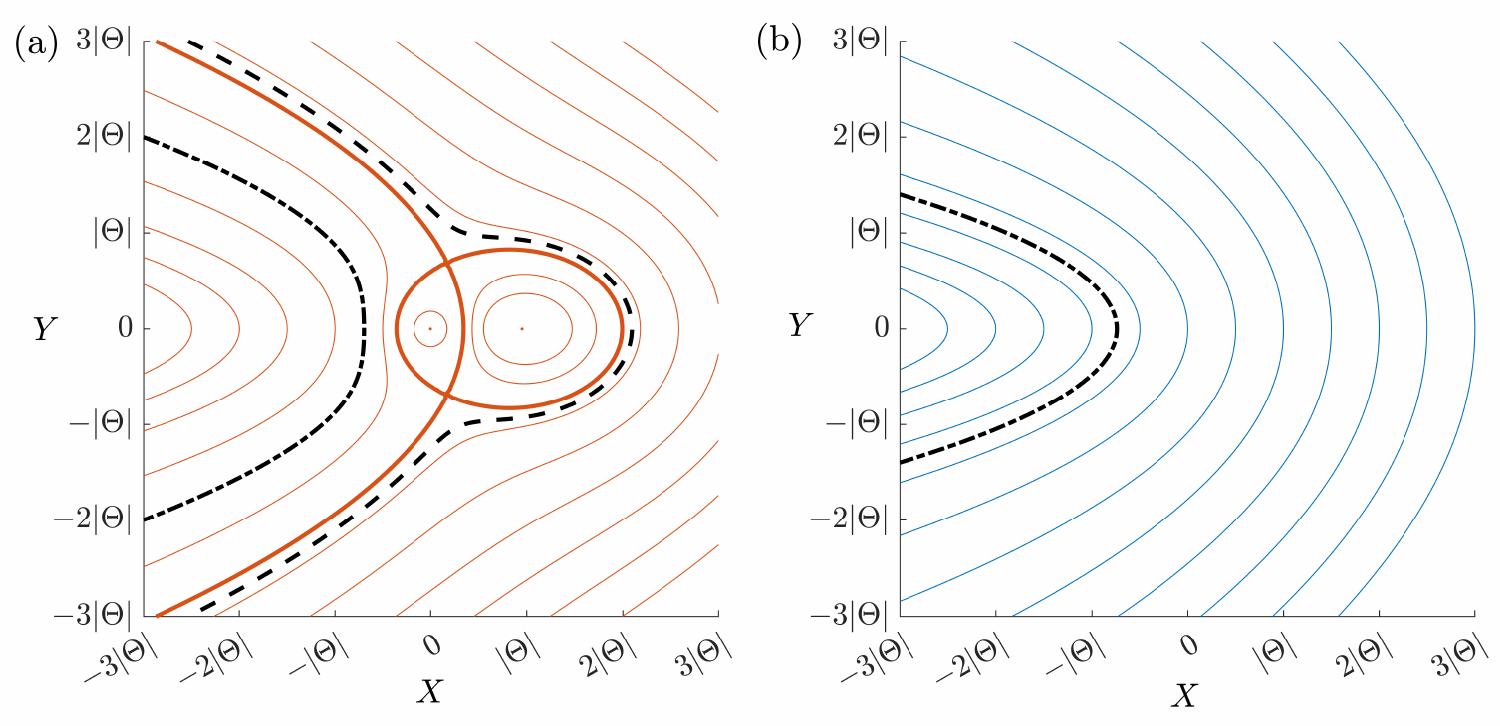}
    \caption{The phase planes for $\Gamma = 0.4 \in \qty(0 , \frac{\sqrt{3}}{2})$. \textbf{(a)} The case $\Theta < 0$ with singular points $\SingPGamma$ (left) and $\SingMGamma$ (right). \textbf{(b)} The case $\Theta > 0$, which has no fixed points or singular points.}    
    \label{fig:gammaless0.8}
\end{figure}

Fig.~\ref{fig:gammaless0.8} shows phase planes for $\Gamma<\frac{\sqrt{3}}{2}$.  The phase plane for $\Theta <0$ is equivalent to that in Fig.~\ref{fig:gammainsqrt1}. However, the $\Theta>0$ phase plane has changed significantly. At $\Gamma = \frac{\sqrt{3}}{2}$, the equilibria $\EqGamma$ and $\EqMinus$ collide and annihilate in a saddle-node bifurcation, so the phase plane contains no equilibria or singularities. All the orbits for $\Theta>0$ are of the (non-extended) direct scattering type.

\subsection{Explaining the scattering for $\Gamma \neq 1 $}
\label{sec:unequal_scattering}
The setup of the three-vortex scattering phenomenon in the generalized system remains as shown in Figure~\ref{fig:setup}, except that $\Gamma_2 = \Gamma \neq 1$ and the two points forming the dipole are separated by a distance $\frac{d}{\Gamma}$ with positions $\r_1 = \pair{-L}{\rho+ \frac{\Gamma d}{2}}$, $\r_2 = \pair{0}{-d}$, and  $\r_3 = \pair{-L}{\rho- \frac{\Gamma d}{2}}$. We will again take $d=1$. The generalized Hamiltonian and angular momentum in the new coordinates, $H \rightarrow  \log{(\Gamma )}$ as $ L \to +\infty,$ and,
\begin{equation*}
       \Theta = \Gamma (1+ 2 \rho ). 
\end{equation*}

\begin{figure}[htbp] 
   \centering
   \includegraphics[width=.45\textwidth]{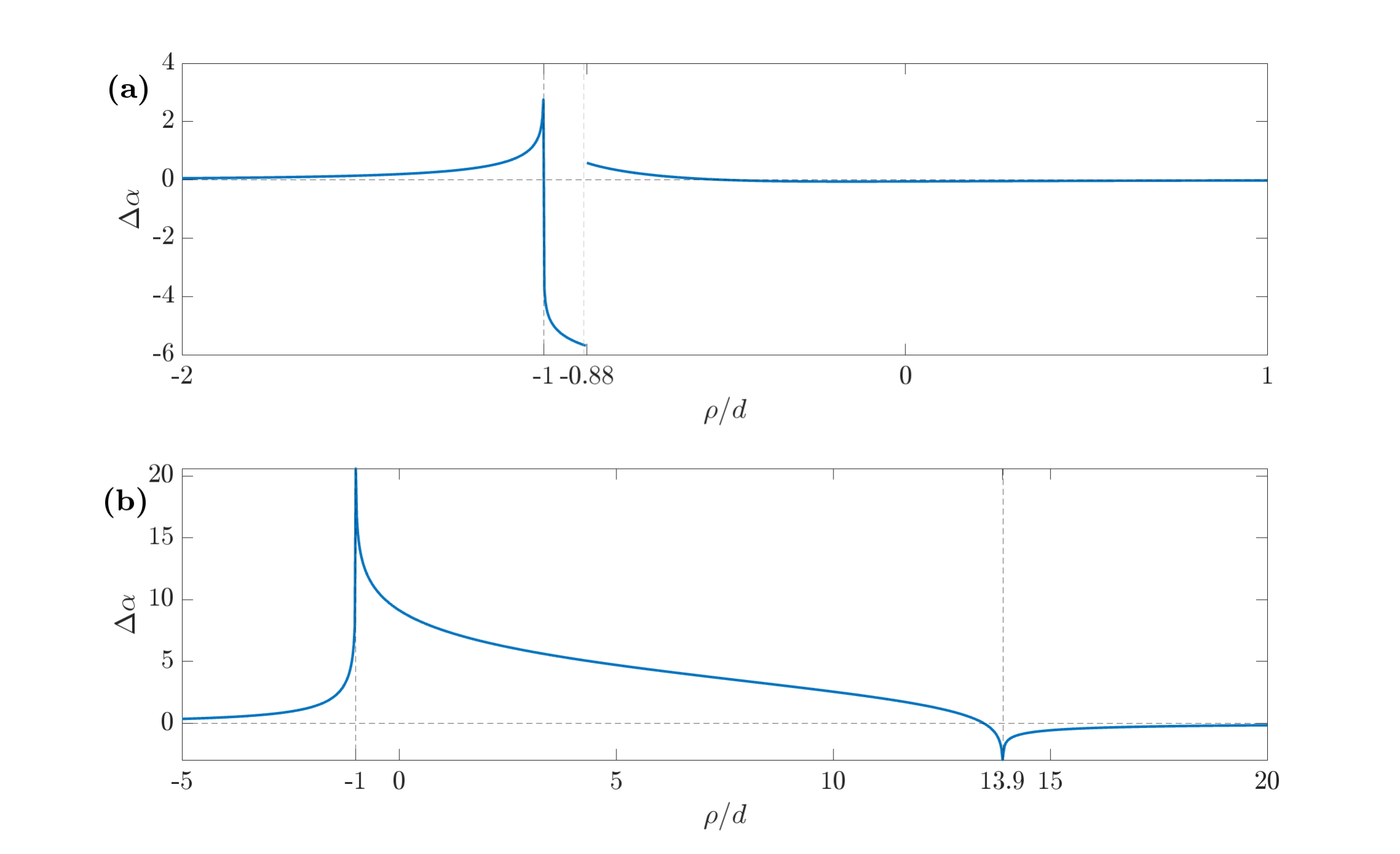} 
\caption{\textbf{(a)} The scattering angle as a function of $\rho$ for $\Gamma=0.4$. \textbf{(b)}~The case $\Gamma=1.7$.}
   \label{fig:theta_vs_rho_generalized}
\end{figure}

\begin{figure}[htbp] 
   \centering
   \includegraphics[width=0.45\columnwidth]{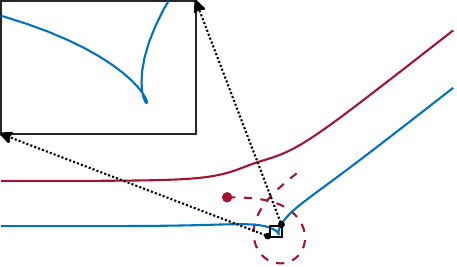} 
   \includegraphics[width=0.45\columnwidth]{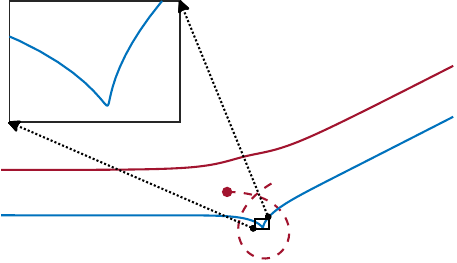} 
   \caption{\label{fig:disappearingloop}\textbf{(a)} The vortex trajectories with $\Gamma=0.4$ and $\rho=-0.9$. \textbf{(b)} The trajectories with $\rho=-0.85$. The insets show a small loop on the trajectory of vortex 3 in the left image that has disappeared in the right image, explaining the $2\pi$ jump in the dependence of the scattering angle shown in Fig.~\ref{fig:theta_vs_rho_generalized}.}
   
\end{figure}

We follow the process described in Section~\ref{sec:scattering} to calculate the critical energy.
For $\Theta <0,$ that is, for $\rho < - \frac{1}{2}$, the critical energy level remains the energy of the equilibria $\EqTriPM$. This again leads to the value $\rhocritminus = -1$ given by Eq.~\eqref{rho_critical_minus}.
For $\Theta > 0$ and $\Gamma>\frac{\sqrt{3}}{2}$, the critical energy is again that of the collinear equilibrium $\EqMinus$
\begin{equation*}
 \rhocritplus = 
   \frac{(\Gamma +1)^2 \left(\Gamma-1\right) \left(\frac{(B+1) (\Gamma +1)^2 \left(\Gamma ^2-1\right)}{-2 A \Gamma ^2+B \left(\Gamma ^4-1\right)-\left(1-\Gamma ^2\right)^2}\right)^{\Gamma }}
   {2 \left(A \Gamma +B \right)} 
    -\frac{1}{2} ,
\end{equation*}
where 
\begin{equation*}
A  = 1-2\Gamma^2 + \sqrt{ 4\Gamma^2 -3} \qand
B  =  \sqrt{\qty(\Gamma^2 -1)^2 + \Gamma^2 A^2}.
\end{equation*}
This matches the value $\frac{7}{2}$ given by~\eqref{rho_critical_plus} for $\Gamma=1$.
Since there is no hyperbolic fixed point when $\Theta>0$ and $\Gamma<\frac{\sqrt{3}}{2}$, there should be only one singular point in the scattering diagram. Two such diagrams are shown in Fig.~\ref{fig:theta_vs_rho_generalized}, demonstrating the disappearance of the second singularity for small $\Gamma$. For $\Gamma=0.4$, the curve jumps by $2\pi$ at $\rho \approx -0.88$. This is explained by the disappearance of a loop in the path of vortex 3; see Fig.~\ref{fig:disappearingloop}.

\section{\label{sec:conclusion} Conclusion}

In this paper, we have introduced a coordinate system for the three-vortex system that, in contrast with previously used reduction methods, avoids introducing artificial singularities into the equations of motion by preserving the topology of the dynamics. They maintain the problem's Hamiltonian structure by introducing Nambu brackets. These coordinates simplify phase-space reasoning and shed new insight into the scattering between a vortex dipole and an isolated vortex. The singular dynamics in trilinear coordinates are equivalent to projecting the dynamics described here into the plane $Y=0$, with the singular curve $\partial \DomainPhysical$ equivalent to the symmetry line $Y=0$ of the spherical or hyperboloidal phase surface.

This reduction should help analyze additional problems in point-vortex dynamics. We mention several such problems.  First, the trilinear coordinate system has been applied to related systems of point vortices, and the reduction developed in this paper should simplify their analysis. A simple example is quasigeostrophic vortices in which a Bessel function replaces the logarithmic potential, but the dynamics is essentially equivalent\cite{Yim:2022}. More interestingly, the four-vortex problem in the integrable case where the total circulation and linear impulse both vanish, has been reduced to trilinear coordinates \cite{Aref:1999}, and here the dynamics become more complicated.

Second, classifying all the changes to the dynamics as the circulations change is surprisingly complicated. Many papers get partway to this goal. Aref first attempted this in the 1979 paper introducing the trilinear coordinate system\cite{Aref.1979}, but the singularity of collinear relative equilibria in these coordinates hampered this effort. Conte classified the bifurcations of the relative equilibria and performed a partial stability analysis\cite{conte2015exact} using a reduction that is very difficult to interpret. Tavantzis and Ting made another study using Synge's trilinear formulation\cite{Tavantzis.1988}. Aref, citing his difficulty in following this analysis, returned to the problem in 2009\cite{Aref.2009c5n}. That approach finds the bifurcations of relative equilibria and their stability but does not describe the dynamics beyond this. The analysis of other phenomena such as the self-similar collapse of the vortex triple has required yet other coordinate systems~\cite{Krishnamurthy:2018,Badin:2018}.  By contrast, the coordinate system introduced here describes the global dynamics in the simplest form possible while yielding equations that can be analyzed using standard methods, even near collapsing states and collinear relative equilibria.

Finally, we mention the motion of four vortices. It is well known that the interaction of two dipoles leads to chaotic scattering\cite{Eckhardt.19884bj,Price.1993,Tophoj.2008}, but the analysis in previous results is somewhat cursory and makes few quantitative predictions. The motion is non-integrable, so Nambu bracket reductions do not apply. However, Ref.~\onlinecite{Price.1993} demonstrates a chaotic scattering process consisting of a sequence of three-vortex interactions in which the fourth vortex remains far from the three strongly interacting vortices during each interaction. Thus, our analysis of the three-vortex problem will serve as the leading-order part of an asymptotic analysis of the problem in this limit. 

\begin{acknowledgments}
The authors gratefully acknowledge support from the NSF under DMS–2206016. We thank Emad Masroor for carefully reading this manuscript and making insightful suggestions for its improvement.
\end{acknowledgments}

\appendix

\section{Detailed calculation of scattering angles}
\label{sec:appendix}

In this appendix, we calculate the change in angle $\Delta \alpha$ on trajectories with initial conditions as $t \to -\infty$ given in Fig.~\ref{fig:setup}. The result is equivalent to one calculated in the supplementary material to~\cite{Lydon.2022}. We include it for completeness and to highlight the connection with the phase planes of Fig.~\ref{fig:threephaseplanes}. 

\begin{figure*}[htbp] 
   \centering
   \includegraphics[width=0.85\textwidth]{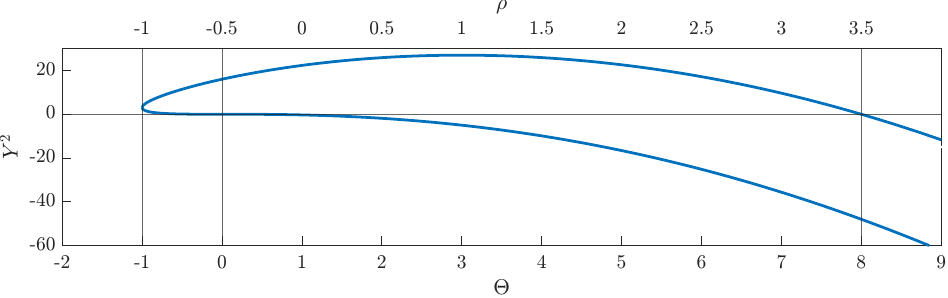} 
   \caption{ The solutions to $p_4(Y^2,\Theta)=0$, with the transitions between the factored form in Eq.~\eqref{p4cases} marked be vertical lines. }
   \label{fig:YsqrTheta}
\end{figure*}

To obtain an explicit integral form, we divide $\dv{\alpha}{t}$ from Eq.~\eqref{alphadot} by $\dv{Y}{t}$, given by Eq.~\eqref{Ydot}, yielding $\dv{\alpha}{Y}$. We remove the dependence on $X$ and $Z$ using the conservation laws~\eqref{HXYZ} and~\eqref{ThetaZXYminus}, and then replace $H$ by its value given the initial condition in Fig.~\ref{fig:setup}. We will use $\Theta$ instead of $\rho$ as the parameter in what follows because it gives somewhat simpler formulas and can use Eq.~\eqref{Theta_rho} to rewrite this in terms of the parameter $\rho$ defining the initial conditions. Integrating this, we find
\begin{equation} \label{psi2Integral}
\Delta \alpha = \int_{\Ymin}^{\infty} \frac{ -8 \Theta^2 \dd Y }{(Y^2+\Theta^2) \sqrt{p_4(Y^2;\Theta)}} \\
  + \int_{\Ymin}^{\infty}\frac{ 8 (\Theta^2 - 8 \Theta) \dd Y }{(Y^2+\Theta^2 - 8 \Theta) \sqrt{p_4(Y^2;\Theta)}} ,
\end{equation}
where
$$
p_4(Y^2;\Theta) = Y^4+2 \left(\Theta^2-4 \Theta -8\right) Y^2 +(\Theta -8) \Theta^3.
$$
These are \emph{complete elliptic integrals}~\cite{Byrd:1971}. To place them in standard form, we must first factor $p_4(Y^2;\Theta)$. We plot its zero locus in Fig.~\ref{fig:YsqrTheta} as a function of $\Theta$ and $Y^2$. From this image, it is clear that $p_4$ can be factored as follows
\begin{widetext}
\begin{equation} \label{p4cases}
p_4(Y^2,\Theta) = \begin{cases}
 (Y^2 - (a+i b)^2)(Y^2 - (a-i b)^2 ),& a>0, b>0, \text{ if } \Theta < -1; \\
 (Y^2 - a^2) (Y^2-b^2), & a>b>0, \text{ if } -1 < \Theta < 0 ; \\
 (Y^2 - a^2) (Y^2+b^2), & a>0, b>0,  \text{ if } 0 < \Theta < 8; \\
 (Y^2 + a^2) (Y^2+b^2), & a>b>0, \text{ if }8 < \Theta.
\end{cases}
\end{equation}
\end{widetext}
The first two cases correspond to the left phase plane of Fig.~\ref{fig:threephaseplanes}, the last two to the right phase plane; the first and last cases correspond to direct scattering, and the second and third to exchange scattering. The lower limit of integration is $\Ymin=0$ in the first and fourth cases, while in the second and third $\Ymin=a$. 
Both integrals in Eq.~\eqref{psi2Integral} can be evaluated with the help of references such as Gradshteyn/Ryzhik and Byrd/Friedman\cite{Gradshteyn:2014,Byrd:1971}. It is quite possible that these expressions can be simplified further. For example, Lydon derived formulas in which $\alpha$ is the sum of one complete elliptic integral of the first kind and one of the third kind.

In the four regions, the constants evaluate to the following
$$
\binom{a^2}{b^2} = \begin{cases}
\frac{1}{2}\begin{pmatrix} \sqrt{\Theta -8} \Theta^{3/2}-\Theta^2+4 \Theta +8 \\
      \sqrt{\Theta -8} \Theta^{3/2}+\Theta^2-4 \Theta -8 \end{pmatrix}
      & \text{ if } \Theta < -1 ; \\
\begin{pmatrix}-\Theta^2+4 \Theta +8 \sqrt{\Theta +1}+8 \\
      -\Theta^2+4 \Theta -8 \sqrt{\Theta +1}+8 \end{pmatrix}
       & \text{ if } -1 < \Theta < 0  ; \\
\begin{pmatrix}-\Theta^2+4 \Theta +8 \sqrt{\Theta +1}+8 \\
      \Theta^2-4 \Theta +8 \sqrt{\Theta +1}-8 \end{pmatrix}
       & \text{ if } 0 < \Theta < 8; \\
\begin{pmatrix}\Theta^2-4 \Theta +8 \sqrt{\Theta +1}-8\\
      \Theta^2-4 \Theta -8 \sqrt{\Theta +1}-8 \end{pmatrix}
      & \text{ if }8 < \Theta.
\end{cases}
$$
In each of the four $\rho$ intervals, the scattering angle can be written as a linear combination of complete elliptic integrals of the first kind 
\begin{align*}
K(m) & = \int_0^1 \frac{\dd x}{\sqrt{\left(1-x^2\right)\left(1-m x^2\right)}} \\
&= \int_0^{\frac{\pi}{2}} \frac{\dd \theta}{\sqrt{1-m \sin^2{\theta}}}.
\end{align*}
and the third kind
\begin{align*}
\Pi(n,m) &= \int_0^1 \frac{\dd x}{\left(1-n x^2\right) \sqrt{\left(1-x^2\right)\left(1-m x^2\right)}} \\
&= \int_0^{\frac{\pi}{2}} \frac{\dd \theta}{\left(1-n \sin^2{\theta}\right)\sqrt{1-m \sin^2{\theta}}}.
\end{align*}

The convention is to define these functions for $0<m<1$, though they are analytic for all $m$ except for a branch cut from $m=1$ to $m=\infty$. 

We report the values found in each of the cases.

\section*{\textbf{Direct scattering with} $\rho<-1$}
Here $\Theta<-1$, and
$$
\Delta \alpha = \frac{64 \sqrt[4]{\Theta -8}\left(-K(m) + \Pi(n,m)\right)}{\sqrt[4]{\Theta } \left(\Theta -8 +\sqrt{\Theta^2 - 8 \Theta }\right) \left(\Theta +\sqrt{\Theta^2 - 8 \Theta }\right)}$$
with
$$
m = \frac{1}{2}+ \frac{4-\Theta }{2 \sqrt{\Theta^2 - 8 \Theta} } \qand
n = \frac{1}{2}-\frac{\Theta^2-4 \Theta -8}{2 \Theta\sqrt{\Theta^2 - 8 \Theta} }.
$$

\section*{\textbf{Exchange scattering with} $-1<\rho<-\frac{1}{2}$}
In this case $-1<\Theta<0$, and
$$ 
 \Delta\alpha = \frac{4 \Theta K(m) + 8\sqrt{1+\Theta} \Bigl(\Pi(n_1,m) -\Pi(n_2,m)\Bigr)}{\sqrt{-\Theta^2+4 \Theta +8 \sqrt{\Theta +1}+8}},   
$$
where
\begin{align*}
m& = \frac{8 +4 \Theta-\Theta^2  -8 \sqrt{\Theta +1}}{8 +4 \Theta -\Theta^2 +8 \sqrt{\Theta +1}}, \\
n_1 &= \frac{\Theta -2 +2 \sqrt{\Theta +1}}{\Theta +2-2 \sqrt{\Theta +1}},\\
 \qand
n_2 &=\frac{\Theta +2 -2 \sqrt{\Theta +1}}{\Theta + 2+2 \sqrt{\Theta +1}}.
\end{align*}
\section*{\textbf{The borderline case} $\rho=-\frac{1}{2}$}
This is the case $\Theta=0$ discussed in Fig.~\ref{fig:rho_minus_half}. Vortex 2 travels along a straight line with no deflection, so the scattering angle is $\alpha=0$.

\section*{\textbf{Exchange scattering with} $-\frac{1}{2}<\rho<\frac{7}{2}$}
Here $0<\Theta<8$, and
$$
\Delta\alpha = \frac{-\Theta^2+4 \Theta +8 \sqrt{\Theta +1}+8}{2 \sqrt[4]{\Theta +1} \left(\Theta +2 \sqrt{\Theta +1}+2\right)} \Bigl(\Pi(n_1,m)-\Pi(n_2,m)\Bigr),$$
where

\begin{align*}
m &= \frac{1}{2} + \frac{\Theta^2-4\Theta-8}{16\sqrt{1-\Theta}}; \\
n_1 &= \frac{2 - \Theta - 2 \sqrt{1+\Theta}}{4}; \\
\qand n_2 &= \frac{2 + \Theta - 2 \sqrt{1+\Theta}}{4}.  
\end{align*}

\section*{\textbf{Direct scattering with} $\frac{7}{2}<\rho$} In this last case, $\Theta>8$ and
$$
\Delta \alpha = c_K K(m) + c_{\Pi,1}\Pi_1(n_1,m) + c_{\Pi,2}\Pi(n_2,m),
$$
where
\begin{widetext}
\begin{gather*}
m = \frac{16 \sqrt{\Theta +1}}{\Theta^2-4 \Theta -8+8 \sqrt{\Theta +1}},\
n_1 = -\frac{4}{\Theta +2 \sqrt{\Theta +1}-2}, \
n_2 = \frac{4 \left(\Theta +2 \sqrt{\Theta +1}+2\right)}{\Theta^2},\\
c_K = -\frac{4 \Theta }{\sqrt{\Theta^2-4 \Theta +8 \sqrt{\Theta +1}-8}}, \\
c_{\Pi,1} = \frac{-2 \Theta^3+4 \Theta^2+64 \Theta +64-4 \sqrt{\Theta +1} \left(\Theta^2-8 \Theta -16\right)}{\sqrt{(\Theta -8) \Theta^3 \left((\Theta -4) \Theta -8 \left(\sqrt{\Theta +1}+1\right)\right)}},\\
\qand c_{\Pi,2} = \frac{-2 \Theta^3+12 \Theta^2-32 \Theta -64 +4 \left(\Theta^2-16\right) \sqrt{\Theta +1}}{\sqrt{(\Theta -8) \Theta^3 \left((\Theta -4) \Theta -8 \left(\sqrt{\Theta +1}+1\right)\right)}}.
\end{gather*}
\end{widetext}

\bibliography{references}

\end{document}